\documentclass[preprint,12pt]{elsarticle}

\usepackage{amssymb}
\usepackage{amsmath}
\usepackage{amsthm}
\usepackage{thmtools}

\usepackage{amsfonts,mathtools,bm,enumitem,cases}
\usepackage[colorlinks=true]{hyperref}
\usepackage[capitalise,nameinlink,noabbrev]{cleveref}
\usepackage{booktabs,xspace}
\usepackage{multirow}

\usepackage{siunitx}
\usepackage[scale=0.8]{geometry}
\usepackage{xcolor}

\newcommand{\pd}{\mathcal{\partial}}

\newcommand{\imh}{{i - \frac 1 2}}
\newcommand{\iph}{{i + \frac 1 2}}
\newcommand{\imhm}{{i - \frac 1 2, -}}
\newcommand{\imhp}{{i - \frac 1 2, +}}
\newcommand{\iphm}{{i + \frac 1 2, -}}
\newcommand{\iphp}{{i + \frac 1 2, +}}
\newcommand{\imhinline}{{i - 1/2}}
\newcommand{\iphinline}{{i + 1/2}}
\newcommand{\imhminline}{{i - 1/2, -}}
\newcommand{\imhpinline}{{i - 1/2, +}}
\newcommand{\imhpminline}{{i - 1/2, \pm}}
\newcommand{\iphminline}{{i + 1/2, -}}
\newcommand{\iphpinline}{{i + 1/2, +}}
\newcommand{\iphpminline}{{i + 1/2, \pm}}
\newcommand{\dx}{{\Delta x}}
\newcommand{\dt}{{\Delta t}}
\newcommand{\jump}[1]{\left[#1\right]}

\newcommand{\ARS}{A-RS\xspace}

\newcommand{\HLL}{\text{HLL}}
\newcommand{\eqvar}{\text{eq}}
\newcommand{\QQ}{\Xi}
\newcommand{\EE}{\xi}

\newcommand{\FWBscheme}{\texttt{FWB1}\xspace}
\newcommand{\FWBschemeTwo}{\texttt{FWB2}\xspace}
\newcommand{\HLLscheme}{\texttt{HLL}\xspace}
\newcommand{\ideal}{\texttt{ideal}\xspace}
\newcommand{\vdW}{\texttt{vdW}\xspace}
\newcommand{\RedKwo}{\texttt{R-K}\xspace}
\newcommand{\PenRob}{\texttt{P-R}\xspace}
\newcommand{\water}{\texttt{H\textsubscript{2}O}\xspace}
\newcommand{\methane}{\texttt{CH\textsubscript{4}}\xspace}
\newcommand{\reldiff}{\eta}

\DeclareFontFamily{U}{mathx}{}
\DeclareFontShape{U}{mathx}{m}{n}{<-> mathx10}{}
\DeclareSymbolFont{mathx}{U}{mathx}{m}{n}
\DeclareMathAccent{\widehat}{0}{mathx}{"70}
\DeclareMathAccent{\widecheck}{0}{mathx}{"71}

\makeatletter
\def\widebreve{\mathpalette\wide@breve}
\def\wide@breve#1#2{\sbox\z@{$#1#2$}%
    \mathop{\vbox{\m@th\ialign{##\crcr
                \kern0.08em\brevefill#1{0.8\wd\z@}\crcr\noalign{\nointerlineskip}%
                $\hss#1#2\hss$\crcr}}}\limits\hspace*{-0.175em}}
\def\brevefill#1#2{$\m@th\sbox\tw@{$#1($}%
    \hss\resizebox{#2}{\wd\tw@}{\rotatebox[origin=c]{90}{\upshape(}}\hss$}
\makeatletter

\newtheorem{theorem}{Theorem}[section]
\newtheorem{lemma}[theorem]{Lemma}

\theoremstyle{definition}
\newtheorem{definition}[theorem]{Definition}

\theoremstyle{remark}
\newtheorem{remark}[theorem]{Remark}

\usepackage{tikz}
\usetikzlibrary{arrows.meta, calc, quotes, tikzmark, shapes.misc, positioning, arrows, decorations.pathreplacing}
\tikzset{cross/.style={cross out, draw=black, minimum size=2*(#1-\pgflinewidth), inner sep=0pt, outer sep=0pt},
    cross/.default={1pt}}
\usepackage{pgfplots}
\pgfplotsset{compat=1.18}
\usepgfplotslibrary{colormaps}
\usetikzlibrary{external}
\tikzsetexternalprefix{../figures}

\pgfplotsset{
	colormap={rainbow desaturated}{
			rgb=(0.278431, 0.278431, 0.858824),
			rgb=(0, 0, 0.360784  ),
			rgb=(0, 1, 1  ),
			rgb=(0, 0.501961, 0  ),
			rgb=(1, 1 , 0  ),
			rgb=(1, 0.380392 , 0  ),
			rgb=(0.419608, 0  , 0 ),
			rgb=(0.878431, 0.301961  , 0.301961 ),
		},
}

\pgfplotsset{
	/pgfplots/surf shading/precision=pdf,
}

\makeatletter
\newsavebox\myboxA
\newsavebox\myboxB
\newlength\mylenA

\newcommand*\xoverline[2][0.75]{%
    \sbox{\myboxA}{$\m@th#2$}%
    \setbox\myboxB\null
    \ht\myboxB=\ht\myboxA%
    \dp\myboxB=\dp\myboxA%
    \wd\myboxB=#1\wd\myboxA
    \sbox\myboxB{$\m@th\overline{\copy\myboxB}$}
    \setlength\mylenA{\the\wd\myboxA}
    \addtolength\mylenA{-\the\wd\myboxB}%
    \ifdim\wd\myboxB<\wd\myboxA%
    \rlap{\hskip 0.5\mylenA\usebox\myboxB}{\usebox\myboxA}%
    \else
    \hskip -0.5\mylenA\rlap{\usebox\myboxA}{\hskip 0.5\mylenA\usebox\myboxB}%
    \fi}

\makeatother

\numberwithin{equation}{section}

\definecolor{review1}{RGB}{117,112,179}
\definecolor{review2}{RGB}{217,95,2}
\definecolor{review3}{RGB}{27,158,119}

\journal{Computers \& Fluids}

\begin{document}
\tikzexternaldisable

\begin{frontmatter}



    \title{Towards a fully well-balanced and entropy-stable scheme for the Euler equations with gravity: General equations of state}


    \author[1]{Victor Michel-Dansac}
    \author[1]{Andrea Thomann}
    \affiliation[1]{
        organization={Université de Strasbourg, CNRS, Inria, IRMA},
        postcode={F-67000},
        city={Strasbourg},
        country={France}
    }

    \begin{abstract}
        The present work concerns the derivation of a fully well-balanced Godunov-type finite volume scheme for the Euler equations with a gravitational potential based on an approximate Riemann solver {in a one-dimensional framework}. It is an extension to general equations of states of the entropy-stable and fully well-balanced scheme for ideal gases recently forwarded in Berthon et al.~\cite{Berthon2025}.
        {A second-order extension preserving the properties of the first-order scheme is given.}
        {The scheme is provably entropy-stable and positivity-preserving for all thermodynamic variables.}
        Numerical test cases illustrate the performance {and entropy stability} of the new scheme, using six different equations of state as examples, four analytic and two tabulated ones.
    \end{abstract}




\end{frontmatter}

\section{Introduction}
\label{sec:intro}

The Euler equations with gravity are the core of many models in fluid dynamics, e.g. in astrophysics or meteorology.
Therein, the hydrodynamic evolution of atmospheres is often described by the full or barotropic Euler equations.
Since these atmospheres are usually a stable system, they are often in an equilibrium configuration, and processes like convection can be viewed as perturbations of the underlying steady state.
Depending on the event that is observed, these perturbations can vary by several orders of magnitude.
For the numerical resolution of those flows, perturbations whose amplitude is smaller than the background error of the applied numerical scheme are quite challenging to properly approximate, even with high order schemes.
In these cases, the mesh resolution has to be refined to reduce the background error and make the perturbations visible in the numerical simulation.
This leads to a computational overhead, which could be avoided were the numerical scheme able to resolve the underlying equilibrium state with high accuracy, or even at machine precision.

    {Such schemes are known as} well-balanced schemes, a name which arose from the pioneering work \cite{Cargo1994} of Cargo \& LeRoux.
They were the first to construct a scheme for the Euler equations with gravitational source terms which is capable of preserving exactly a discrete form of a hydrostatic equilibrium.
Further, in the work of Bale \& LeVeque \cite{LeVeque1999} within the quasi-steady wave-propagation framework, the source term was numerically included in the Riemann problem, a technique which was extended from shallow-water equations to the Euler equations with an ideal gas equation of state.
Since then, a lot of research has been devoted to the development of well-balanced schemes, more particularly in the context of the shallow water equations.
For instance, \cite{Audusse2004,Pares2004,XinZhaShu2010,Berthon2012c,Chen2017,CabCasMorMun2023} is a non-exhaustive list of contributions to the topic.
For the Euler equations, much work has been dedicated to numerically preserving the special class of stationary hydrostatic atmospheres \cite{GabCasDum2018a,Li2018,BerKapChaKli2021,XinShu2012,Kaeppeli2014};
see also Käppeli's detailed review~\cite{Kaeppeli2022} for an overview of techniques and concepts.

More recently, so-called moving steady states got more and more attention,
first in the context of the shallow water equations \cite{Mantri2024,Castro2007,NoeXinShu2007,MicBerClaFou2016},
but also for general systems \cite{GomezBueno2021,CasPar2020,BerChaKli2021,BerBulFouMbaMic2022,FraMicNav2023}.
For the Euler equations, strategies that were effective for hydrostatic equilibria are often adapted to achieve well-balanced solutions for some moving equilibria. For example, in \cite{KliPupSem2019}, a method was developed to preserve a constant, non-zero velocity in a single grid-aligned direction.

Besides the well-balanced property, a key property for robustness is entropy stability.
For the Euler equations, on the analytical level, an accompanying inequality, the entropy inequality, can be derived.
So-called admissible entropy weak solutions of the Euler equations satisfy this relation, which is important for excluding numerical solutions which are not thermodynamically compatible \cite{Tadmor1986, Tadmor2003}.
To ensure that the numerical solution is the physical admissible entropy solution, the numerical method has to fulfil the discrete analogue of this entropy inequality.
Such schemes are then called entropy-stable.
Their construction is quite challenging and a lot of work was dedicated to deriving entropy-stable schemes for the homogeneous Euler equations
    {\cite{HarLaxLee1983,Chalons2008e,FjoMisTad2012,RayChaFjoMis2016,BerKli2021,Mar2023}
    }
and also entropy-stable well-balanced schemes which preserve hydrostatic equilibria \cite{Desveaux2016,Thomann2019}.
Recently, an entropy-stable fully well-balanced scheme for the Euler equations {under the ideal gas assumption} was constructed in \cite{Berthon2025} which preserves both {isentropic} hydrostatic and moving equilibria.
It belongs to the family of so-called structure-preserving numerical schemes, i.e., schemes where physical properties of the model are preserved within the discrete numerical solution.
Besides fulfilling the entropy inequality, it guarantees the positivity of thermodynamic quantities, namely pressure, temperature, etc.

To link these thermodynamic quantities and close the equations, an additional relation, the so-called equation of state (EOS) is needed.
A standard assumption that can be found in many approaches in the literature is to consider an ideal gas law.
It assumes that the modelled fluid is a gas, whose particles are point-like and do not interact with each other which simplifies the model in contrast to considering a general EOS.
However, this assumption is quite restrictive and limits the applicability of the numerical scheme.
Consequently, for real applications, other EOSs have to be considered.
In \cite{Kaeppeli2014}, a general EOS was considered when well-balancing barotropic hydrostatic equilibria.
This numerical scheme was then extended in \cite{GrosheintzLaval2020} to steady adiabatic flows with grid-aligned streamlines.
Also, in \cite{Varma2019}, hydrostatic equilibria with classical van-der-Waals gases and van-der-Waals gases with radiation pressure were considered.
However, these are only special examples and their extension to general EOS is, at times, quite technical and challenging.

Therefore, this work is dedicated to the construction of a fully well-balanced scheme for general equations of state that is also entropy-stable in the region of thermodynamic compatibility, i.e., in regions where the weak solution of the Euler equations fulfills an accompanying entropy inequality.
    {In the current work we consider the one-dimensional framework only. Firstly, the construction of approximate Riemann solvers is only rigorous in one dimension and secondly, in multiple dimensions, the equilibrium states include a divergence free constraint on the momentum, which requires truly multi-dimensional solvers which can preserve involution constraints on machine precision.
        To our knowledge, this has not been achieved yet for the Euler equations with gravity.}

As in the ideal gas case, the positivity of thermodynamic variables follows from the entropy stability and the provable positivity of the density.
To our knowledge, this is the first time that the entropy stability was rigorously proven and applied to the case of fully well-balanced schemes for general equations of states.
The scheme can be applied on analytically available equations of states or on tabulated ones, which increases its applicability.

The resulting scheme is quite simple to implement, since it amounts to a mere modification of intermediate states in a Riemann solver and source term discretization.
Moreover, despite the EOS being potentially nonlinear,
we highlight that the scheme itself does not feature any nonlinear iteration,
except when computing the EOS itself.

The paper is organized as follows.
The next section is dedicated to the description of the Euler equations with gravity and equilibrium solutions.
This section also contains the definition of the equations of state used in the numerical simulations, namely four with analytical expressions based on a cubic EOS and two tabulated ones.
The numerical scheme is then derived in \cref{sec:WBscheme}.
First, we describe the general framework of the Godunov-type finite volume scheme, and the conditions that have to be fulfilled by the associated approximate Riemann solver are recalled from \cite{Berthon2025} and adjusted to the case of general EOS.
Subsequently, the intermediate states of the Riemann solver are derived, as well as a compatible discretization of the gravitational source term.
In \cref{sec:NumRes} the accuracy and performance of the numerical scheme with respect to the theoretical properties of the numerical scheme are assessed for the above-mentioned EOSs.
Finally, \cref{sec:Conclusion} concludes this work.

\section{The Euler equations with gravity}
\label{sec:model}

We consider, in a one-dimensional setting,
the compressible Euler equations with a gravitational source term,
{see e.g.~\cite{LeVeque2002} for a well-known textbook on the subject.}
They are governed by
\begin{equation}
    \label{eq:EulerG}
    \begin{dcases}
        \pd_t \rho  + \pd_x (\rho u) = 0 ,                             \\
        \pd_t (\rho u)  + \pd_x (\rho u^2 + p) = - \rho \pd_x \varphi, \\
        \pd_t E  + \pd_x ((E + p) u) = - \rho u \pd_x \varphi.
    \end{dcases}
\end{equation}
Therein $\rho > 0$ denotes the density,
$u \in \mathbb{R}$ the velocity
and $E > 0$ the total energy density, which is given by the sum of the internal energy density $\rho e$ and the kinetic energy, that is
\begin{equation}
    \label{eq:TotalEnergy}
    E = \rho e(\tau,s) + \frac{1}{2}\rho u^2.
\end{equation}
{To close the system, we consider a pressure law
$p(\tau, e): \mathbb{R}^+ \times \mathbb{R}^+ \to \mathbb{R}^+$, given in terms of specific volume $\tau = 1/\rho$
and specific internal energy $e$, which will be defined by different equations of state (EOSs) detailed later on.}
The function $\varphi : \mathbb{R} \to \mathbb{R}$ is a given time-independent
continuous gravitational potential.

Moreover, in accordance with the second law of thermodynamics, we assume the existence of a specific entropy $s(\tau,e): \mathbb{R}^+ \times \mathbb{R}^+ \to \mathbb{R}$,
obeying, for some temperature $T(\tau,e) > 0$, the Gibbs relations
\begin{equation}\label{eq:DerivEntropy}
    \frac{\pd s}{\pd \tau} (\tau,e) = -\frac{p(\tau,e)}{T(\tau,e)} < 0
    \text{\qquad and \qquad}
    \frac{\pd s}{\pd e} (\tau,e) = -\frac{1}{T(\tau,e)}<0.
\end{equation}
{Note that this defines the negative of the quantity commonly known as specific physical or thermodynamic entropy and is thus convex \cite{Galgani1968,CallenBook}.}
According to \eqref{eq:DerivEntropy}, for arbitrary $\tau > 0$, the function $e \mapsto s(\tau,e)$ is strictly decreasing.
Therefore, it is injective (and even bijective),
and we can define the inverse function $e(\tau,s)$, used in the definition of the total energy \eqref{eq:TotalEnergy}.
Other quantities of interest include the specific enthalpy $h$ and the specific total enthalpy $H$, respectively defined by
\begin{equation}
    \label{eq:total_enthalpy}
    H = h + \varphi \text{, \quad with \quad} h = \frac{E + p}{\rho}.
\end{equation}

To shorten notation, we write system~\eqref{eq:EulerG}
in the following compact form
\begin{equation*}
    \label{eq:Euler_compact_form}
    \pd_t W + \pd_xF(W) = S(W),
\end{equation*}
where the state vector $W$, flux function $F$ and source term $S$ are respectively defined by
\begin{equation*}
    W =
    \begin{pmatrix}
        \rho \\ \rho u \\ E
    \end{pmatrix}, \quad
    F(W) =
    \begin{pmatrix}
        \rho u \\ \rho u^2 + p \\ u (E + p)
    \end{pmatrix}, \text{\quad and \quad}
    S(W) =
    \begin{pmatrix}
        0 \\ - \rho \pd_x \varphi \\ - \rho u \pd_x \varphi
    \end{pmatrix}.
\end{equation*}
Equipped with this notation, most thermodynamical quantities can be expressed
in terms of $W$ (and $\varphi$ if needed).
For instance, the specific total enthalpy $H$ can be seen as a function $H(W, \varphi)$,
and the specific entropy $s$ as a function $s(W)$.
In \eqref{eq:EulerG}, we require the density $\rho$ and the specific internal energy $e$ to be positive, which by the pressure law induces the positivity of the pressure~$p$.
Thus, the state vector $W$ belongs to the set of admissible states $\Omega$,
defined by
\begin{equation}
    \label{eq:admissible_states}
    \Omega = \left\{
    W \in \mathbb{R}^3 \text{\; such that \;}
    \rho > 0 \text{\ and \ } {e} > 0
    \right\}.
\end{equation}
To reflect the influence of the gravitational source term on the wave structure of the Euler equations, we augment system \eqref{eq:EulerG} by adding the (trivial) equation $\pd_t \varphi = 0$.
The resulting system is hyperbolic, with eigenvalues
\begin{equation*}
    \lambda_{\pm} = u \pm c, \quad
    \lambda_u = u, \quad
    \lambda_0 = 0,
\end{equation*}
where $c$ denotes the speed of sound,
given by
\begin{equation*}
    c =
    \sqrt{\left(\frac{\partial p(\rho, s)}{\partial \rho}\right)_s},
\end{equation*}
where the subscript $s$ denotes the derivative taken at constant entropy.
This wave structure is represented in the left panel of \cref{fig:Riemann_solver}.
Note that, compared to the homogeneous system,
which is obtained by setting the gravitational source terms in \eqref{eq:EulerG} to zero,
there is {an additional} zero eigenvalue which is associated to
the gravitational potential.
In contrast to the homogeneous Euler equations, this leads to a non-ordered wave structure,
since $\lambda_{\pm}$ and $\lambda_u$ can be positive or negative, depending on the flow.


\subsection{Equilibrium solutions}
\label{sec:moving_equilibria}

Time-invariant solutions of \eqref{eq:EulerG} with a non-zero velocity $u \neq 0$,
so-called moving equilibria, are governed by the following system
\begin{subnumcases}{\label{eq:EulerG_Stationary}}
    \label{eq:EulerG_Stationary.rho}
    \pd_x (\rho u) = 0, \\
    \label{eq:EulerG_Stationary.mom}
    \pd_x (\rho u^2 + p) = - \rho \, \pd_x \varphi, \\
    \label{eq:EulerG_Stationary.energy}
    \pd_x ((E + p) u) = - \rho u \, \pd_x \varphi.
\end{subnumcases}
Assuming smooth enough steady solutions fulfilling \eqref{eq:EulerG_Stationary}, we obtain from \eqref{eq:EulerG_Stationary.rho} a constant momentum $q_0 \coloneqq \rho u \neq 0$.
Substituting and dividing by $q_0$
in \eqref{eq:EulerG_Stationary.energy}, we find that
\begin{equation}
    \label{eq:def_constant_enthalpy}
    \frac{E + p}{\rho} + \varphi \eqqcolon H_0,
\end{equation}
where $H_0$ denotes a constant specific total enthalpy.
Substituting $\varphi$ from \eqref{eq:def_constant_enthalpy}
into \eqref{eq:EulerG_Stationary.mom}, we obtain
\begin{equation}
    \label{eq:EulerG_Stationary.mom2}
    \pd_x \! \left( \frac {q_0^2} \rho + p \right)
    =
    \rho \, \pd_x \! \left( \frac{E + p}{\rho} \right).
\end{equation}
Arguing the definition \eqref{eq:TotalEnergy} of $E$, we find
\begin{equation*}
    \frac{E + p}{\rho} = e + \frac{1}{2} \frac{q_0^2}{\rho^2} + \frac{p}{\rho},
\end{equation*}
and so \eqref{eq:EulerG_Stationary.mom2} reduces to
\begin{equation*}
    - \frac {q_0^2} {\rho^2} \pd_x \rho
    +
    \pd_x p
    =
    \rho \, \pd_x e
    -
    \rho \frac {q_0^2} {\rho^3} \pd_x \rho
    +
    \pd_x p
    -
    \frac{p}{\rho} \pd_x \rho.
\end{equation*}
Eliminating and rearranging terms simplifies the above expression to
\begin{equation*}
    - \frac{p}{\rho} \pd_x \rho + \rho \pd_x e = 0.
\end{equation*}
Rewriting this relation in terms of specific volume $\tau$, we obtain
\begin{equation}
    \label{eq:steady_relation_dx_tau_dx_e}
    p \pd_x \tau + \pd_x e = 0.
\end{equation}
However, according to the expressions \eqref{eq:DerivEntropy}
of the entropy derivatives, note that
\begin{equation*}
    \pd_x s
    =
    \frac{\pd s}{\pd \tau} \pd_x \tau
    +
    \frac{\pd s}{\pd e} \pd_x e
    =
    - \frac{p}{T} \pd_x \tau - \frac{1}{T} \pd_x e.
\end{equation*}
Plugging \eqref{eq:steady_relation_dx_tau_dx_e} into the above expression, we obtain
\begin{equation}
    \pd_x s = 0.
\end{equation}

Hence, smooth moving steady solutions for a general equation of state, with $u \neq 0$,
are necessarily isentropic and are characterized by the constant triplet $(q_0,H_0,s_0)$ given by
\begin{equation}
    \label{eq:equilibrium}
    \rho u =: q_0,
    \qquad
    \frac{E + p}{\rho} + \varphi =: H_0,
    \qquad
    s =: s_0,
\end{equation}
where $s_0$ is a constant specific entropy.
In practice, starting from a given triplet $(q_0,H_0,s_0)$ and gravitational potential $\varphi$,
the steady solution in state variables $W$ is obtained
by solving the nonlinear system \eqref{eq:equilibrium} with Newton's method.

A special case is given by time-invariant solutions of \eqref{eq:EulerG} with zero velocity $u=0$, the so-called hydrostatic equilibria.
They fulfil the ordinary differential equation
\begin{equation}
    \partial_x p = - \rho \partial_x \varphi.
\end{equation}
Since this problem is ill-posed, additional assumptions about the dependence between pressure and density are required.
Intensively studied examples are so-called isothermal or isentropic hydrostatic equilibria with constant temperature or constant entropy, respectively.
The reader is referred to the review \cite{Kaeppeli2022} for a non-exhaustive list of numerical methods achieving well-balancing of hydrostatic equilibria using different numerical techniques.
    {Note that, in our case,
        the moving equilibrium~\eqref{eq:equilibrium}
        has been obtained by assuming that $u \neq 0$,
        in order to divide by $u$
        in~\eqref{eq:EulerG_Stationary.energy}.
        It turns out that taking $u = 0$ in~\eqref{eq:equilibrium}
        yields the isentropic hydrostatic equilibrium,
        which can also be obtained
        from~\eqref{eq:EulerG_Stationary},
        assuming a vanishing velocity at constant entropy.}

\subsection{Equation of state}
\label{sec:equation_of_state}

In this work, we will apply various equations of state (EOS), including both analytical expressions and tabulated forms.

First, we consider a general cubic EOS, with the analytical form
\begin{equation}
    \label{eq:EOS_cubic}
    \begin{aligned}
        p(\tau,T)  & = \frac{R T}{\tau - b} - \frac{a(T)}{(\tau - b r_1) (\tau - b r_2)}
        \text{\qquad and \qquad}                                                           \\
        e(\tau, T) & = c_v (\tau,T) T + \frac{a(T) - T a^\prime(T)}{b} \mathfrak{u}(\tau),
    \end{aligned}
\end{equation}
where $R$ denotes the gas constant, $b$ the covolume, $r_1,r_2 \in \mathbb{R}$ are parameters,
and $\mathfrak{u}(\tau)$ is a function specific to each EOS.
The function $a(T)$ defines an attraction term, and $c_v$ is the heat capacity at constant volume. It depends on the temperature and the specific volume as follows
\begin{equation}
    c_v(\tau,T) = c_v^0 - \frac{T a^{\prime\prime}(T)}{b} \mathfrak{u}(\tau),
\end{equation}
where $c_v^0$ is a constant.
The gas constant $R$ is connected to the ratio of specific heats $\gamma$ by the relation $R = (\gamma - 1) c_v$.
Depending on the definition of the free parameters and functions, one obtains different known equations of state, such as the ideal gas law, van-der-Waals (\vdW) EOS \cite{VanDerWaals1873,Toro2009}, Redlich-Kwong (\RedKwo) EOS \cite{Redlich1949} or the Peng-Robinson (\PenRob) EOS \cite{Peng1976}.
The details of the free parameters and functions are given in \cref{tab:EOS}.
Moreover, the corresponding {specific entropy} definition reads
\begin{equation}\label{eq:def.entropy}
    s(\tau,T) = -\left(s_0 - \frac{a^{\prime}(T)}{b} \mathfrak{u}(\tau) + R \log\left(\tau - b\right) + c_v^0 \log\left(T\right)\right).
\end{equation}
Note that, in the numerical scheme, the {specific entropy} is required to be given in terms of~$\tau$ and~$e$.
It can be recovered analytically using \eqref{eq:def.entropy} and the temperature formulas given in \cref{tab:EOS_temperature}.
However, $e(\tau,s)$ is also required by the numerical scheme.
Unfortunately, recovering the internal energy in terms of $\tau$ and $s$ involves solving a nonlinear equation, which is done numerically by employing Newton's method (except in the \ideal and \vdW cases, where a closed form is available).
    {Typically, Newton's method takes between $4$ and~$8$ iterations to reach a precision of $10^{-14}$.}

\begin{table}[!t]
    \renewcommand{\arraystretch}{1.5}
    \centering
    \begin{tabular}{ccccc}
        \toprule
        EOS     & $r_1$           & $r_2$           & $a(T)$                                                      & $\mathfrak{u}(\tau)$                                                \\
        \cmidrule(lr){1-5}
        \ideal  & 0               & 0               & 0                                                           & 0                                                                   \\
        \vdW    & 0               & 0               & $a_0$                                                       & $- \dfrac{b}{\tau}$                                                 \\
        \RedKwo & 0               & $-1$            & $\dfrac{a_0}{\sqrt{T}}$                                     & $\log\Big(\dfrac{\tau}{\tau +b}\Big)$                               \\
        \PenRob & $-1 - \sqrt{2}$ & $-1 + \sqrt{2}$ & $a_0\left(1 + \kappa\left(1 - \sqrt{T/T_0}\right)\right)^2$ & $\dfrac{1}{r_1 - r_2} \log\Big(\dfrac{\tau-br_1}{\tau -b r_2}\Big)$ \\
        \bottomrule
    \end{tabular}
    \caption{Parameters and defining functions for some cubic EOS, in particular the ideal gas  (\ideal), van-der-Waals (\vdW), Redlich-Kwong (\RedKwo) and Peng-Robinson (\PenRob) EOS. Therein, $a_0$, $\kappa$, $c_v^0$ and $b$ denote constants which will be specified when setting up the numerical test cases.}
    \label{tab:EOS}
\end{table}

\begin{table}[!tb]
    \renewcommand{\arraystretch}{2}
    \centering
    \begin{tabular}{crcl}
        \toprule
        \ideal  & $T(\tau,e)$ & = & $\dfrac{e}{c_v^0}$                                                                                                                                                                 \\
        \cmidrule(lr){1-4}
        \vdW    & $T(\tau,e)$ & = & $\dfrac{e - a_0 \mathfrak{u}(\tau)}{c_v^0}$                                                                                                                                        \\
        \cmidrule(lr){1-4}
        \RedKwo & $T(\tau,e)$ & = & $\dfrac{{{\left(\QQ\left( \tau ,e\right) +\sqrt[3]{4} e\right) }^{2}}}{3 \, \sqrt[3]{4} \, c_v^0 \, \QQ\left( \tau ,e\right) }$,                                                   \\
                &             &   & with \quad $\QQ\left( \tau ,e\right) = {{\left( \sqrt{27 {{\EE\left( \tau \right) }^{2}} {{c_v^0}}-4 {{e}^{3}}}- \EE\left( \tau \right)  \sqrt{27 c_v^0}\right) }^{\frac{2}{3}}} $ \\
                &             &   & where \quad $\EE(\tau) = \dfrac 3 2 a_0 \mathfrak{u}(\tau)$                                                                                                                        \\
        \cmidrule(lr){1-4}
        \PenRob & $T(\tau,e)$ & = & $\dfrac{1}{4 T_0 (c_v^0)^2}
            \left(
            \sqrt{
                    {{\EE\left( \tau \right) }^{2}} {{\kappa}^{2}} +
                    4 T_0 c_v^0 \Big(
                    (\kappa + 1) \EE(\tau) + e
                    \Big)
                }-\EE\left( \tau \right) \kappa
        \right)^2$,                                                                                                                                                                                                    \\
                &             &   & with \quad $\EE(\tau) = a_0 (\kappa+1) \mathfrak{u}(\tau) \vphantom{\dfrac 1 2}$                                                                                                   \\
        \bottomrule
    \end{tabular}
    \caption{Temperatures in terms of $\tau$ and $e$ for some cubic EOS: cases of the ideal gas (\ideal), van-der-Waals (\vdW), Redlich-Kwong (\RedKwo) and Peng-Robinson (\PenRob) EOS.}
    \label{tab:EOS_temperature}
\end{table}

Further, we consider two tabulated EOS from the \texttt{CoolProp} library \cite{CoolProp}, which contains the thermodynamic description of pure and pseudo-pure fluids and mixtures.
As examples in the numerical experiments,
we choose the IAPWS-95 formulation for water~\cite{WagPru2002},
as well as the EOS for methane~\cite{SetWag1991}.

\subsection{Entropy inequality}
\label{sec:entropy_inequality}

We begin this section by defining the
\emph{mathematical entropy} for a hyperbolic system
such as~\eqref{eq:EulerG}.
This definition contrasts with the
thermodynamically motivated \emph{specific entropy} $s$,
defined in \eqref{eq:DerivEntropy}.
A convex real valued function $U(W)$ is called a mathematical entropy for \eqref{eq:EulerG} if there exists a real valued function $G(W)$, called the entropy flux, which satisfies
\begin{equation*}
    \nabla_W U(W)^\intercal \nabla_W F(W) = \nabla_W G(W)^\intercal,
\end{equation*}
where $F$ denotes the flux defined in \eqref{eq:Euler_compact_form}.
For the Euler system \eqref{eq:EulerG} we have the family of mathematical entropies given by $U(W) := \rho \eta(s(W))$  with entropy fluxes $G(W) = \rho u \eta(s(W))$.
Therein, $\eta$ is a function such that $W \mapsto U(W)$ is convex with respect to the state variables $W$.
A sufficient condition on $\eta$ is that $\eta'(s) \geq 0$ and $\eta''(s) \geq 0$.
This has been proven for a general equation of state under the usual thermodynamic assumptions in \cite[Section 1.4]{Bouchut2004}.
Therefore, admissible entropy weak solutions to system \eqref{eq:EulerG} satisfy the entropy inequality (in the sense of distributions)
\begin{equation}
    \label{eq:entropy_ineq_id}
    \pd_t (\rho \eta(s))
    +
    \pd_x (\rho\eta(s) u) \leqslant 0.
\end{equation}

In the following, the objective is to approximate the system \eqref{eq:EulerG} with an arbitrary EOS.
To that end, we derive a Godunov-type finite volume (FV) scheme based on an approximate Riemann solver (\ARS), which is able to preserve {isentropic} moving {and hydrostatic} equilibria up to machine precision,
generates admissible solutions in the sense of \eqref{eq:admissible_states},
and fulfills all entropy inequalities~\eqref{eq:entropy_ineq_id} {discretely}.

\section{The numerical scheme}
\label{sec:WBscheme}

We begin with a brief overview of the principle behind Godunov-type finite volume (FV) schemes, relying on approximate Riemann solvers.
As is standard in the finite volume framework, we divide the computational domain $I\subset \mathbb{R}$ into non-overlapping cells $C_i = (x_\imh, x_\iph)$ with center $x_i$.
    {The size of cell $C_i$ is then given by \smash{{$\dx = x_\iph - x_\imh$}}, which is chosen for simplicity uniformly for all cells composing the grid. }
Then, at a given time $t^n$, the volume average of the solution on cell $C_i$ is defined as
\begin{equation}
    \label{eq:wn}
    W_i^{n}=\frac{1}{\dx}\int_{x_\imh}^{x_\iph}
    W(x,t^{n})dx.
\end{equation}
The time line is discretized in time steps of variable size $\Delta t$. We focus here on one generic time step $[t^n , t^{n+1} ]$ where $t^{n+1} = t^n + \Delta t$.
Since the states $W^{n+1}$ rely only on data at $t^n$, i.e. an explicit scheme, the time step has to be restricted by a CFL (Courant-Friedrichs-Lewy) condition,
see~\cite{HarLaxLee1983}, which relates $\Delta t$ to the cell size $\Delta x$ and ensures that the information is propagated to a maximum of one cell.

In the next section, we develop an approximate Riemann solver (\ARS) for a general EOS, {similar to the \ARS based on the ideal gas EOS in} \cite{Berthon2025}.

\subsection{The approximate Riemann solver}
\label{sec:Riemann_solver}

In the following, we define an \ARS, denoted by~$\widetilde W$,
whose role is to approximate the solution of the Riemann problems
occurring at each cell interface $x_\iphinline$ between cells $C_i$ and $C_{i+1}$.
Since this solution is self-similar,
the \ARS~depends on the variable $x/t$,
and on the two neighboring states $W_i^n$ and $W_{i+1}^n$
located on the left and right of the interface, respectively.
That is to say, the approximate Riemann solution
at interface $x_\iphinline$ is given by
$\smash{\widetilde W((x - x_\iphinline) / t; W_i^n, W_{i+1}^n)}$.

In this work, we choose an approximate Riemann solution consisting of two intermediate states,
as shown in the right panel of \cref{fig:Riemann_solver}.
{Note that the middle wave corresponds here to a wave with velocity zero, unlike other well-known approximate Riemann solvers, like the HLLC solver from~\cite{Toro1994}, where the middle wave is an approximation of the contact wave with velocity $u$. As a consequence, the present structure may lead to diffused contact waves.
The expression of the approximate Riemann solution}, for two given left and right states $W_L$ and $W_R$, is as follows
\begin{equation}
    \label{eq:ARS}
    \widetilde{W}\left( \frac x t; W_L, W_R \right) =
    \begin{dcases}
        W_L   & \text{if } x < -\lambda t,     \\
        W_L^* & \text{if } -\lambda t < x < 0, \\
        W_R^* & \text{if } 0 < x < \lambda t,  \\
        W_R   & \text{if } x > \lambda t,      \\
    \end{dcases}
\end{equation}
where the components of the intermediate state vectors $W_L^*$ and $W_R^*$ read
\begin{equation}
    \label{eq:W_star}
    W_L^* =
    \begin{pmatrix}
        \rho_L^* \\ q_L^* \\ E_L^*
    \end{pmatrix}
    \text{\qquad and \qquad}
    W_R^* =
    \begin{pmatrix}
        \rho_R^* \\ q_R^* \\ E_R^*
    \end{pmatrix}.
\end{equation}
{Note that due to the gravitational source terms, the intermediate states depend not only on the states $W_L, W_R$ but also on the gravitational potential averaged in the sense of \eqref{eq:wn} respectively denoted by $\varphi_L, \varphi_R$.}
{To define the scheme based on the \ARS \eqref{eq:ARS},
it is useful to introduce its cell-wise juxtaposition
at time $t^{n+1}$.
For all $i \in \mathbb{Z}$,
for all $x \in (x_i, x_{i+1})$, we define
\begin{equation*}
    W_\Delta^{n+1}(x)
    =
    \widetilde W \left(
    \frac{x - x_\iph}{\dt};
    W_i^n, W_{i+1}^n
    \right).
\end{equation*}
Then, the scheme is obtained by integrating
this juxtaposition over the cell
$C_i = (x_\imhinline, x_\iphinline)$, to get
\begin{equation}
    \label{eq:scheme_from_juxtaposition}
    W_i^{n+1}
    =
    \frac 1 \dx \int_{x_\imh}^{x_\iph}
    W_\Delta^{n+1}(x) \, dx.
\end{equation}
Thanks to the specific structure of the \ARS \eqref{eq:ARS},
following e.g. \cite{MicBerClaFou2016,Berthon2025},
\eqref{eq:scheme_from_juxtaposition} rewrites
\begin{equation}
    \label{eq:scheme_from_ARS}
    W_i^{n+1} = W_i^n - \lambda \frac \dt \dx \Big[
        W_L^* \left( W_i^{n}, W_{i+1}^n, \varphi_i, \varphi_{i+1} \right)
        -
        W_R^* \left( W_{i-1}^{n}, W_i^n, \varphi_{i-1}, \varphi_i \right)
        \Big],
\end{equation}
which can be rewritten into the standard conservation form.}

\begin{figure}[tb]
    \centering
    \begin{tikzpicture}[scale=0.8]
        \pgfmathsetmacro{\dx}{-7}
        \pgfmathsetmacro{\eps}{0.2}

        \draw[thick, -latex] (-2.5+\dx,0) -- (2.5+\dx,0) node[right]{$x$};
        \draw[thick, -latex] (0+\dx,-0.25) -- (0+\dx,3.25) node[above]{$t$};

        \draw[] (\dx,0) -- (\dx+0.5,2.5) node[above]{$u$};
        \draw (\dx,0) -- (\dx+1.75,2.5) node[above]{$u+c$};
        \draw (\dx,0) -- (\dx-0.75,2.5) node[above]{$u-c$};

        \draw[thick, -latex] (-2.5,0) -- (2.5,0) node[right]{$x$};
        \draw[thick, -latex] (0,-0.25) -- (0,3.25) node[above]{$t$};
        \draw (0,0) -- (-1.75-\eps, 2.5) node[above]{$-\lambda$};
        \draw (0,0) -- (1.75+\eps, 2.5) node[above]{$+\lambda$};

        \node at (-1.6, 0.9) {$W_L$};
        \node at (-0.6, 2) {$W_L^*$};
        \node at (0.6, 2) {$W_R^*$};
        \node at (1.6, 0.9) {$W_R$};
    \end{tikzpicture}
    \caption{%
        Left panel: One possible wave configuration,
        with $u > 0$,
        of the Euler equations with gravity \eqref{eq:EulerG}.
        Right panel: Wave structure of the approximate Riemann solver.%
    }
    \label{fig:Riemann_solver}
\end{figure}

In this work, we impose a symmetry on the waves defining the \ARS around \(\lambda_0 = 0\) {using a single speed $\lambda$}.
This is a simplifying assumption and in principle, two different values e.g. {$\lambda_L$ and $\lambda_R$}, could be chosen to
capture the asymmetric wave structure of the Euler equations with respect to \(\lambda_0\),
provided that the acoustic waves \(\lambda_\pm\) remain within the \ARS.
In practice, we set
\begin{equation}
    \label{eq:def_lambda}
    \lambda = \Lambda \, \max \Big(|u_L| + c(W_L), |u_R| + c(W_R) \Big), \quad \Lambda \geq 1.
\end{equation}
Note that this expression
does not necessarily satisfy
$\lambda \geq |\lambda_\pm|$ for the choice of  $\Lambda = 1$,
see~\cite{TorMueSiv2020}, however it is frequently applied and is sufficient to ensure stability and accuracy for the numerical scheme.
The role of the parameter~$\Lambda$ is therefore
to scale the wave speeds to increase stability at the cost of slightly increasing the artificial diffusion of the solver.
Its value (typically $1$) will be given in the numerical result section for each numerical experiment.

Equipped with $\lambda$, we are able to define the CFL restriction on $\dt$.
Indeed, let \smash{$\lambda_\iph$} be the approximate wave speed at the interface \smash{$x_\iph$},
computed by applying \eqref{eq:def_lambda} to the states $W_i^n$ and~$W_{i+1}^n$.
For large enough $\Lambda$, the inequality
\smash{$\lambda_\iph \geq |\lambda_\pm|$} holds.
Under this assumption,
the scheme is stable as long as
the time step satisfies
\smash{$\dt \, \max_i \lambda_\iph \leqslant \dx$}.
However, the consistency proof, see the theorems
from~\cite{HarLaxLee1983}, requires the
more restrictive condition
\begin{equation}
    \label{eq:CFL_condition}
    \dt \leqslant \frac 1 2 \frac{\dx}{\max_i \lambda_\iph},
\end{equation}
which guarantees that interactions between Riemann solutions are avoided.

\subsection{Consistency conditions}
\label{sec:Riemann_solver_consistency}

At this level, we have introduced six unknown quantities in \eqref{eq:W_star},
namely two intermediate states each for the density, momentum and energy stored in $W_L^*$ and $W_R^*$.

In this section, we revisit conditions on the intermediate
states~$W_L^*$ and $W_R^*$ to ensure the consistency of the FV scheme with solutions of the Euler system.
These conditions are classical, see \cite{HarLaxLee1983},
and were derived in detail in~\cite{Berthon2025}.
The main idea is to ensure that the
integral of the approximate Riemann solution
on the space-time domain $(x_i,x_{i+1})\times [t^n,t^{n+1}]$
is equal to that of the exact Riemann solution.
In the following we motivate and summarize the consistency result from \cite{Berthon2025}.
The wave structure of the \ARS we develop here can be interpreted as the one of the classical HLL Riemann solver from Harten, Lax and van Leer~\cite{HarLaxLee1983} with the additional zero wave associated to the gravitational potential.
Therefore, it is natural to extend the HLL solver to include the effects of gravity while ensuring the well-balancedness and preserving the positivity and entropy stability of the HLL solver.

Thus, to write the intermediate states, we start from the HLL intermediate states, modified by the presence of the source term.
This leads to the following expressions for the intermediate states,
\begin{equation}
    \label{eq:IShatplusdeviation}
    W_L^* = \widehat W - \delta W, \quad W_R^* = \widehat W + \delta W,
\end{equation}
where, due to the Harten-Lax consistency conditions from \cite{HarLaxLee1983}, $\widehat W$ satisfies
\begin{equation}
    \label{eq:expression_w_hat}
    \left\{
    \begin{aligned}
        \widehat \rho & = \rho_\HLL,   \\
        \widehat q    & = q_\HLL
        + \frac {S^q \dx} {2 \lambda}, \\
        \widehat E    & = E_\HLL
        + \frac {S^E \dx} {2 \lambda}. \\
    \end{aligned}
    \right.
\end{equation}
The intermediate state $W_\HLL$ in \eqref{eq:expression_w_hat}
is the intermediate state of the HLL solver, given by
\begin{equation}
    \label{eq:def_w_hll}
    W_\HLL =
    \begin{pmatrix}
        \rho_\HLL \\ q_\HLL \\ E_\HLL
    \end{pmatrix} =
    \frac{W_L + W_R} 2 - \frac{F(W_R) - F(W_L)} {2 \lambda}.
\end{equation}
Moreover, we have introduced two terms $S^q$ and $S^E$,
which are consistent approximations of the cell averages
of the source terms $S^{(2)}(W) = -\rho \pd_x \varphi$ and
$S^{(3)}(W) = -q \pd_x \varphi$, respectively, and will be defined later.

In \cite{Berthon2025} it was shown that
choosing the intermediate states
\eqref{eq:IShatplusdeviation}
ensures the consistency of the resulting \ARS \eqref{eq:ARS}.
Hence, the scheme \eqref{eq:scheme_from_ARS} and its conservative counterpart are consistent as well.
We refer the interested reader to e.g. \cite{HarLaxLee1983,MicBerClaFou2016} for details.

To complete the description of the numerical scheme, the goal is thus to derive the remaining five unknowns
(namely, $\delta \rho$, $\delta q$, $\delta E$, $S^q$ and $S^E$)
while ensuring positivity, entropy stability and well-balancedness.
These properties will be defined in the next section.

\subsection{Properties and conditions to be satisfied by the \ARS}
\label{sec:definitions}

We begin with the definition of the well-balanced property.
We are interested in developing a so-called
\emph{fully well-balanced} scheme,
which exactly preserves the (moving) steady solutions
described in \cref{sec:moving_equilibria} associated with isentropic equilibria.
These states are characterized by constant momentum $q$,
entropy $s(W)$, and specific total enthalpy $H(W, \varphi)$, given by \eqref{eq:total_enthalpy}.
{We distinguish two concepts, the interface steady solution (ISS) which is a local characterization regarding the \ARS, and well-balancedness, which is a global characterization at the level of the numerical scheme \eqref{eq:scheme_from_ARS}.}

\begin{definition}
    \label{def:ISS}
    A pair $(W_L, W_R)$ of admissible states is said to
    be an \emph{Interface Steady Solution} (ISS) if, and only if,
    \begin{equation*}
        q_L = q_R
        \text{, \quad}
        H(W_L, \varphi_L) = H(W_R, \varphi_R)
        \text{, \quad and \quad}
        s(W_L) = s(W_R).
    \end{equation*}
    Then, the scheme \eqref{eq:scheme_from_ARS} is said to be well-balanced if
    \begin{equation*}
        \forall i \in \mathbb{Z}, \; (W_i^n, W_{i+1}^n) \text{ is an ISS}
        \qquad \implies \qquad
        \forall i \in \mathbb{Z}, \; W_i^{n+1} = W_i^n.
    \end{equation*}
\end{definition}

The following result regarding the connection between ISS and intermediate states
were proven in \cite{Berthon2025}.

\begin{lemma}[Sufficient condition for well-balancedness]
    \label{cor:WB_consequence}
    Let $(W_L, W_R)$ be an ISS.
    A sufficient condition for well-ba\-lan\-ced\-ness is
    \begin{equation}
        \label{eq:WB_ISS}
        W_L^* = W_L, \quad \text{and} \quad W_R^* = W_R.
    \end{equation}
    Since the triplet $(q,H,s)$ is constant for a steady state, this condition implies
    \begin{equation}
        \label{eq:WB_consequence_on_intermediate_states}
        q^\ast = q_L^* = q_R^*
        \text{, \quad}
        H(W_L^*, \varphi_L) = H(W_R^*, \varphi_R), \quad
        s^\ast = s(W_L^*) = s(W_R^*),
    \end{equation}
    and implies that the intermediate states of the \ARS \eqref{eq:IShatplusdeviation} satisfy
    \begin{equation}
        \label{eq:WB_consequence_on_hat_and_delta_w}
        \widehat W = \frac{W_L + W_R} 2
        \text{ \; and \; }
        \delta W = \frac{W_R - W_L} 2.
    \end{equation}
\end{lemma}

Conditions \eqref{eq:WB_consequence_on_intermediate_states}
and \eqref{eq:WB_consequence_on_hat_and_delta_w}
will be crucial in determining the expressions for
\(\delta W\) and \(S^q\),~\(S^E\) in
\cref{sec:determination_of_delta_w,sec:source_term_approximation}, respectively.

Next, we turn to the condition for entropy stability.
    {Note that the definition of entropy stability is a mathematically motivated concept going beyond the thermodynamic compatibility. }

\begin{definition}[Entropy stability]
    \label{def:entropy_stability}
    The scheme \eqref{eq:scheme_from_ARS} under the CFL condition \eqref{eq:CFL_condition} is said to be entropy-stable if for all (mathematical) entropies $\rho\eta(s)$ satisfying \eqref{eq:entropy_ineq_id},
    for all $i\in\mathbb{Z}$,
    the states $W_i^n$ and $W_i^{n+1}$ fulfil the discrete entropy inequality
    \begin{equation}
        \label{discreteentropineq}
        \rho_i^{n+1}\eta(s_i^{n+1}) \leqslant \rho_i^n \eta(s_i^n)
        - \frac{\dt}{\dx} \left((\rho \eta(s) u)_\iph^n -
        (\rho \eta(s) u)_\imh^n\right).
    \end{equation}
\end{definition}
The entropy stability result proven in \cite{Berthon2025}, although given in the context of the ideal gas law, does not make any assumption on the EOS and can thus be applied to the case of general EOSs.
It only relies on the well-known integral entropy consistency, stated in \cite{HarLaxLee1983},
{and is expressed in terms of the
        HLL averages of the entropy and of the entropy flux,
        defined with $s_L = s(W_L)$ and $s_R = s(W_R)$ by
        \begin{align*}
            (\rho s)_\HLL                 & = \frac{1}{2}
            \big(\rho_Ls_L + \rho_Rs_R\big)
            -\frac{1}{2\lambda}\big(\rho_Rs_Ru_R -\rho_Ls_Lu_L\big), \\
            \big(\rho \eta(s)\big)_{\HLL} & = \frac{1}{2}
            \big(\rho_L\eta(s_L) + \rho_R\eta(s_R)\big)
            -\frac{1}{2\lambda}\big(\rho_R\eta(s_R)u_R -\rho_L\eta(s_L)u_L\big).
        \end{align*}}
In \cite{Berthon2025} is proven that the entropy stability for the numerical scheme is obtained if, and only if, the \ARS satisfies
\begin{equation}
    \label{eq:Entropyconsistency}
    \rho_\HLL \eta(s^*) \leqslant (\rho \eta(s))_\HLL.
\end{equation}
Further, it is proven that for all smooth convex functions $\eta$, the choice
\begin{equation}
    \label{eq:condition_on_s_star}
    s^* = \frac{(\rho s)_\HLL}{\rho_\HLL}
\end{equation}
satisfies the entropy consistency condition \eqref{eq:Entropyconsistency} and thus the numerical scheme is entropy stable.
    {Note that the definition of the intermediate entropy \eqref{eq:condition_on_s_star} is only well-defined if $\rho_\HLL$ is positive.}
    {However, as detailed in \cite{Berthon2025}, this holds as soon as $\lambda > \max(|u_L|,|u_R|)$, which in our case is satisfied by the definition of $\lambda$ given in \eqref{eq:def_lambda}.}

The final property we want to preserve is the positivity of the thermodynamic variables, here in terms of density and internal energy. The positivity of e.g. temperature and pressure depends then on the chosen EOS.

\begin{definition}[Positivity preservation]
    \label{def:positivity_preservation}
    Let $W_i^n \in \Omega$ for all $i \in \mathbb{Z}$ and $n \geq 0$.
    The scheme \eqref{eq:scheme_from_ARS} is said to be positivity-preserving if,
    for all $i \in \mathbb{Z}$, $W_i^{n+1} \in \Omega$,
    where $\Omega$ is the set of admissible states defined in \eqref{eq:admissible_states}.
\end{definition}
According to the definition \eqref{eq:scheme_from_juxtaposition}
of $W_i^{n+1}$, the positivity of the intermediate densities
$\rho_L^*$ and~$\rho_R^*$ is a necessary condition
for the positivity of $\rho_i^{n+1}$.
Since the intermediate states for the density are given by $\rho_\HLL$ and a fluctuation $\delta\rho$, their positivity is ensured as soon as the fluctuation does not exceed $\rho_\HLL$ itself.
This is essentially the statement of the next result whose proof can be found in \cite{Berthon2025}. As it does not use any assumption on the EOS it also holds for general EOSs.
\begin{lemma}[Density positivity]
    Let $W_L \in \Omega$ and $W_R \in \Omega$ two given states.
        {Let $\lambda$ be defined by \eqref{eq:def_lambda}, which ensures that $\rho_{\HLL}>0$.}
    Then,
    \begin{equation}
        \label{eq:positivity_condition_on_rho}
        |\delta \rho| < \rho_\HLL
    \end{equation}
    is a sufficient condition for $\rho_L^\ast > 0$ and $\rho_R^\ast >0 $.
\end{lemma}
{The positivity of the thermodynamic variables is
closely related to the entropy stability.
Note that due to the ideal gas assumption in \cite{Berthon2025} which connects directly $\tau, p$ and $s$, therein, the positivity of the pressure was proven which yields the positivity of the internal energy as a consequence.
Since the Gibbs relations directly yield
a relation between $\tau$, $e$ and $s$,
we show here the positivity of the internal energy
instead of the pressure.
This is also motivated by the fact that the transform between internal energy and entropy at a fixed $\tau$ is used in the numerical scheme to define intermediate states of the \ARS.
The proof of the following Lemma \ref{lem:positivity_inte} is along the lines of the positivity proof for the pressure given in \cite{Berthon2025}.
However, since that proof relied on the particular case of the ideal gas law,
we give the details to suit the general EOS case.
\begin{lemma}[Positivity of internal energy]
    \label{lem:positivity_inte}
    Assume that, for all $i \in \mathbb{Z}$
    and for some $n \geq 0$,
    $W_i^n \in \Omega$.
    Further, assume that \eqref{eq:positivity_condition_on_rho}
    is satisfied, and that $\rho_i^{n+1} > 0$.
    Then, $e_i^{n+1} > 0$
    as soon as, for each interface $x_\iphinline$,
    $e_L^* > 0$ and $e_R^* > 0$.
\end{lemma}
\begin{proof}
    Let $i \in \mathbb{Z}$ and $n \geq 0$ be fixed.
    First, recall from \cref{sec:entropy_inequality}
    that the function $W \mapsto U(W) = \rho \eta(s(W))$ is convex.
    Therefore, applying this function to the definition
    \eqref{eq:scheme_from_juxtaposition} of $W_i^{n+1}$
    and using Jensen's inequality, we obtain
    \begin{equation*}
        \label{eq:Jensen_juxtaposition_1}
        \rho_i^{n+1} \eta(s(W_i^{n+1}))
        \leqslant
        \frac 1 \dx \int_{x_\imh}^{x_\iph}
        \rho_\Delta^{n+1}(x) \,
        \eta \bigl( s \bigl(W_\Delta^{n+1}(x) \bigr) \bigr)
        \, dx,
    \end{equation*}
    that is to say
    \begin{align}
        \label{eq:Jensen_juxtaposition_2}
        \begin{split}
            \eta(s(W_i^{n+1}))
            \leqslant &
            \int_{x_\imh}^{x_\iph}
            \eta \bigl( s \bigl(W_\Delta^{n+1}(x) \bigr) \bigr)
            \frac {\rho_\Delta^{n+1}(x) \, dx} {\rho_i^{n+1} \dx} \\
            =         &
            \int_{x_\imh}^{x_\iph}
            \eta \bigl( s \bigl(W_\Delta^{n+1}(x) \bigr) \bigr)
            d \varrho(x),
        \end{split}
    \end{align}
    where we have defined the probability measure
    \begin{equation*}
        d \varrho(x) = \frac {\rho_\Delta^{n+1}(x) \, dx} {\rho_i^{n+1} \dx}.
    \end{equation*}
    Now, for some $\tau > 0$,
    recall that the function $e \mapsto s(\tau, e)$
    is decreasing by Gibbs' relation
    and convex, see e.g. \cite{Galgani1968}
    or \cite[Section 8.1]{CallenBook} for the concavity of $-s$.
    Moreover, as mentioned in \cref{sec:entropy_inequality},
    the function $s \mapsto \eta(s)$ is
    increasing and convex by hypothesis,
    see \cite{Bouchut2004}.
    Therefore, the composition of these two functions
    is convex and decreasing.
    Applying its inverse which is concave and decreasing
    to \eqref{eq:Jensen_juxtaposition_2},
    and using Jensen's inequality again, yields
    \begin{equation*}
        e_i^{n+1}
        \geqslant
        \int_{x_\imh}^{x_\iph}
        e_\Delta^{n+1}(x)
        d \varrho(x).
    \end{equation*}
    Now, at the level of one interface, note that
    $e_\Delta$ is a juxtaposition of the two states
    $e_L$ and $e_R$, and of the intermediate states
    $e_L^*$ and $e_R^*$.
    Both states $e_L$ and $e_R$ are positive by hypothesis,
    and so the positivity of $e_\Delta$ is ensured
    as soon as the intermediate internal energies
    $e_L^*$ and $e_R^*$ are positive.
    The proof is thus concluded.
\end{proof}
}

\subsection{\texorpdfstring{Determination of $\delta W$}{Determination of delta W}}
\label{sec:determination_of_delta_w}

The objective of this section is to utilize the definitions and conditions
introduced in \cref{sec:definitions} to derive an expression
for the three unknown components of \(\delta W\).

We begin with the sufficient conditions for well-balancedness given in condition \eqref{eq:WB_consequence_on_intermediate_states}.
Since the momentum component is constant in a steady state, we have only one intermediate momentum $q^*$.
From \eqref{eq:WB_ISS}, it immediately follows that $q_L = q_R$, which yields $\delta q = 0$.

Analogously, we obtain, for the specific total enthalpy, $H(W_L, \varphi_L) = H(W_R, \varphi_R)$.
Expanding the specific total enthalpy according to its definition \eqref{eq:total_enthalpy},
we obtain a relation on the specific enthalpy $h$ given by
\begin{equation}
    \label{eq:WB_h}
    \frac{E_L + p_L}{\rho_L} + \varphi_L = \frac{E_R + p_R}{\rho_R} + \varphi_R
    \quad \iff \quad
    [\varphi] = - [h].
\end{equation}
From this we obtain an expression for $\delta \rho$ which is independent of the considered EOS,
{and thus merely taking the definition from \cite{Berthon2025} would lead to the desired well-balanced property.
However, the expression for $\delta \rho$ in \cite{Berthon2025}
is not well-defined when $[h] = 0$,
and more specifically when $[\varphi] = [h] = 0$,
since it involves the fraction ${[\varphi]}/{[h]}$.}

{To correct this shortcoming,
we here propose a new expression for $\delta \rho$,
defined in all situations.
This new expression relies on two conditions:
\begin{subequations}
    \begin{align}
        \label{eq:condition_on_delta_rho_1}
        [\varphi] = - [h]
         & \implies \delta \rho = \frac{[\rho]}{2},              \\
        \label{eq:condition_on_delta_rho_2}
        [\varphi] = 0
         & \implies \delta \rho = 0. \vphantom{\frac{[\rho]}{2}}
    \end{align}
\end{subequations}
Condition \eqref{eq:condition_on_delta_rho_1}
is nothing but the well-balancedness requirement,
while condition  \eqref{eq:condition_on_delta_rho_2}
corresponds to recovering the
HLL solver in the absence of a gravitational source term,
which may increase stability thanks to the good
properties of the HLL solver.
However, these conditions are fundamentally incompatible when $[\varphi] = [h] = 0$,
as $\delta \rho$ would both vanish and be non-zero.
Therefore,
we relax condition \eqref{eq:condition_on_delta_rho_2},
to make sure that the scheme remains
well-balanced in all situations.
The new conditions to be satisfied then read
\begin{subequations}
    \begin{align*}
        [\varphi] = - [h]
         & \implies \delta \rho = \frac{[\rho]}{2},              \\
        [\varphi] = 0 \text{ and } [h] \neq 0
         & \implies \delta \rho = 0. \vphantom{\frac{[\rho]}{2}}
    \end{align*}
\end{subequations}
This reduces to making an explicit
choice for the case where $[\varphi] = [h] = 0$.}

The new expression is then based on the following lemma,
introducing a function $\psi$ and its properties.
\begin{lemma}
    \label{lem:helper_function_for_density}
    Let $\psi$ be the function defined by
    \begin{equation*}
        \psi(x, y, \alpha) = \psi_1 \big( \psi_2 (x, y) \big)^\alpha,
    \end{equation*}
    where $\psi_1$ is defined by
    \begin{equation*}
        \psi_1(z) = \cos \left(\frac \pi 2 z\right) \exp \left(-2 z^2 \right),
    \end{equation*}
    and where $\psi_2$ is defined by
    \begin{equation}
        \label{eq:def_psi_2}
        \psi_2(x, y) = \frac{x + y}{\mathcal{M}(\varepsilon_0, \sqrt{x^2+y^2})}.
    \end{equation}
    In \eqref{eq:def_psi_2}, $\varepsilon_0$ is set to $10^{-12}$,
    and $\mathcal{M}$ is a regularized maximum function defined by
    \begin{equation*}
        \mathcal{M}(\varepsilon, z) =
        \begin{dcases}
            \varepsilon       & \text{ if } z < \frac{\varepsilon}{2},              \\
            z                 & \text{ if } z > \frac{3\varepsilon}{2},             \\
            P(\varepsilon, z) & \text{ otherwise}, \vphantom{\frac{\varepsilon}{2}}
        \end{dcases}
    \end{equation*}
    where $P$ is a polynomial defined such that $\mathcal{M}(\varepsilon_0, \cdot) \in \mathcal{C}^2(\mathbb{R}_+^*, \mathbb{R}_+^*)$:
    \begin{equation*}
        P(\varepsilon, z)
        =
        \frac {-z^4} {2 \varepsilon^3}
        + \frac {2 z^3} {\varepsilon^2}
        - \frac {9 z^2} {4 \varepsilon}
        + z
        + \frac {27 \varepsilon} {32}.
    \end{equation*}
    Then, for all $\alpha \in \mathbb{N}$ and $x, y \in \mathbb{R}^2$,
    $\psi$ satisfies the following properties:
    \begin{align}
        \label{eq:psi_i} \tag{$\psi$-i}
        \psi(x, y, \alpha) = 1        & \iff x = - y,                                                            \\
        \label{eq:psi_ii} \tag{$\psi$-ii}
        y \geq 3 \, \varepsilon_0 / 2 & \implies \psi(0, y, \alpha) = 0,                                         \\
        \label{eq:psi_iii} \tag{$\psi$-iii}
        x \geq 3 \, \varepsilon_0 / 2 & \implies \psi(x, 0, \alpha) = 0,                                         \\
        \label{eq:psi_iv} \tag{$\psi$-iv}
        |\psi(x, y, \alpha)|          & \leqslant 1,                                                             \\
        \label{eq:psi_v} \tag{$\psi$-v}
        \psi(\cdot, \cdot, \alpha)    & \in \mathcal{C}^2(\mathbb{R} \times \mathbb{R}, \mathbb{R}),             \\
        \label{eq:psi_vi} \tag{$\psi$-vi}
        y \geq 3 \, \varepsilon_0 / 2 & \implies \psi(x, y, \alpha) \underset{x \to 0}{=} \mathcal{O}(x^\alpha).
    \end{align}
\end{lemma}

\begin{figure}[!ht]
    \centering
    \includegraphics{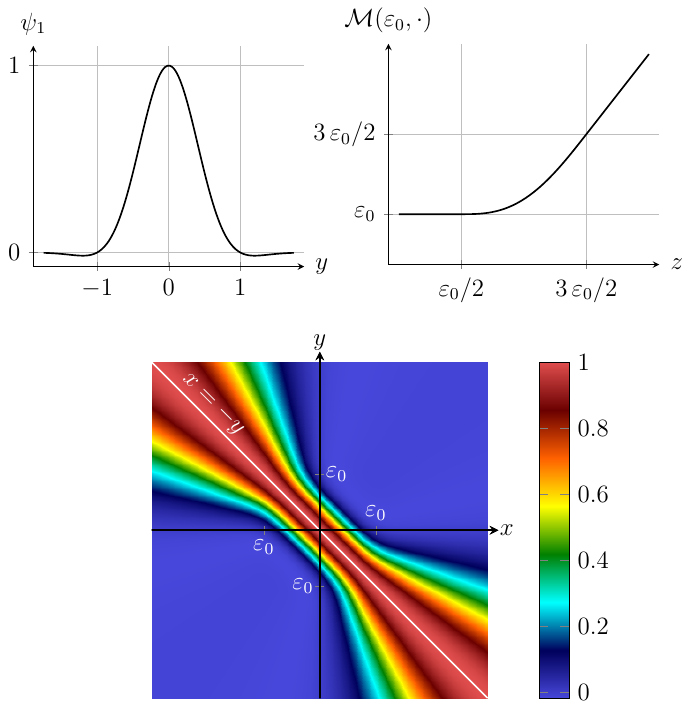}

    \caption{Drawings of the functions $\psi_1$ (top left panel), $\mathcal{M}$ (top right panel), and $\psi$ (bottom panel).}
    \label{fig:psi}
\end{figure}

\begin{proof}
    The properties \eqref{eq:psi_i}--\eqref{eq:psi_iv} are straightforward,
    and come from the properties of $\psi_1$ (proven in \cite{Berthon2025}).
    A sketch of $\psi$ is provided in \cref{fig:psi} to help the reader understand the function.

    To prove \eqref{eq:psi_v}, we rely on the values of $P$ and its derivatives:
    \begin{equation*}
        \begin{aligned}
            P\left(\varepsilon, \frac \varepsilon 2 \right) = \frac \varepsilon 2,
             & \qquad
            P\left(\varepsilon, \frac {3 \varepsilon} 2 \right) = \frac {3 \varepsilon} 2, \\
            \frac{\partial P}{\partial z}
            \left(\varepsilon, \frac \varepsilon 2 \right) = 0,
             & \qquad
            \frac{\partial P}{\partial z}
            \left(\varepsilon, \frac {3 \varepsilon} 2 \right) = 1,                        \\
            \frac{\partial^2 P}{\partial z^2}
            \left(\varepsilon, \frac \varepsilon 2 \right) = 0,
             & \qquad
            \frac{\partial^2 P}{\partial z^2}
            \left(\varepsilon, \frac {3 \varepsilon} 2 \right) = 0.                        \\
        \end{aligned}
    \end{equation*}
    Therefore, $\mathcal{M}$ is a $\mathcal{C}^2$ function, and so is $\psi_2$.
    Since $\psi_1$ is smooth, $\psi$ turns out to be $\mathcal{C}^2$.

    Finally, to prove \eqref{eq:psi_vi}, assume that $y \geq 3 \, \varepsilon_0 / 2$.
    Therefore, for all $x \in \mathbb{R}$,
    \begin{equation*}
        \mathcal{M}(\varepsilon_0, \sqrt{x^2+y^2})
        =
        \sqrt{x^2+y^2}
        =
        |y| + \mathcal{O}(x^2).
    \end{equation*}
    Hence,
    \begin{equation*}
        \psi(x, y, \alpha)
        =
        \psi_1 \left( \frac{x + y}{|y| + \mathcal{O}(x^2)} \right)^\alpha
        =
        \psi_1 \left( \operatorname{sgn}(y) + \mathcal{O}(x) \right)^\alpha.
    \end{equation*}
    However, recall from \cite{Berthon2025} that
    $\psi_1 (\pm 1 + x) = \mathcal{O}(x)$ as $x \to 0$.
    Therefore, we have proven that $\psi(x, y, \alpha) = \mathcal{O}(x^\alpha)$ as $x \to 0$,
    which concludes the proof.
\end{proof}

Equipped with $\psi$, we are now able to state the expression of $\delta \rho$.

\begin{lemma}
    \label{lem:positivity_density}
    Let $W_L \in \Omega$ and $W_R \in \Omega$ two given states,
    and let $W_\HLL$ be given by \eqref{eq:def_w_hll} with {$\lambda$ defined by \eqref{eq:def_lambda}}.
    Further, let $\delta \rho$ be defined by
    \begin{equation*}
        \label{eq:deltarho_definition}
        \delta \rho =
        \frac{[\rho]} 2 \, \psi \big( [\varphi], [h], 1 \big).
    \end{equation*}
    Then, $\delta \rho$ satisfies the following properties:
    \begin{align}
        \label{eq:delta_rho_i} \tag{i}
        [\varphi] = -[h]                                              & \implies \delta \rho = [\rho] / 2, \\
        \label{eq:delta_rho_ii} \tag{ii}
        ([\varphi] = 0 \text{\; and \;} [h] \geq 3\, \varepsilon_0/2) & \implies \delta \rho = 0,          \\
        \label{eq:delta_rho_iii} \tag{iii}
        ([h] = 0 \text{\; and \;} [\varphi] \geq 3\, \varepsilon_0/2) & \implies \delta \rho = 0,          \\
        \label{eq:delta_rho_iv} \tag{iv}
        |\delta\rho|                                                  & \leqslant \rho_{\HLL}.
    \end{align}
\end{lemma}

\begin{proof}
    Properties \eqref{eq:delta_rho_i}--\eqref{eq:delta_rho_iv} are straightforward to prove.
    Indeed, they are directly related to the
    properties \eqref{eq:psi_i}--\eqref{eq:psi_iv} of $\psi$,
    proven in \cref{lem:helper_function_for_density}.
\end{proof}

\begin{remark}
    Property \eqref{eq:delta_rho_i} corresponds to the ISS case,
    and ensures that $\delta \rho$ satisfies
    \eqref{eq:WB_consequence_on_hat_and_delta_w}, while property \eqref{eq:delta_rho_ii} corresponds to the no-gravity case,
    and proves that $\delta \rho$ vanishes,
    leading to $\rho_L^* = \rho_R^* = \rho_\text{HLL}$.
    Finally, property \eqref{eq:delta_rho_iii} ensures that $\delta \rho$ is well-defined
    and vanishes when $[h]$ vanishes.
    This property involved a limit in \cite{Berthon2025},
    whereas the new expression allows it to be directly satisfied.
    Note that, appropriately, both properties \eqref{eq:delta_rho_ii} and \eqref{eq:delta_rho_iii}
    have been relaxed when both $[\varphi]$ and $[h]$ vanish,
    by introducing $\varepsilon_0$.
    Property~\eqref{eq:delta_rho_iv} is the required condition
    \eqref{eq:positivity_condition_on_rho} from
    \cref{lem:positivity_density} for positivity of the density.
\end{remark}

The remaining component $\delta E$ depends on the EOS,
and thus it will be different from the one derived
in \cite{Berthon2025} in the case of the ideal gas law,
and will be part of the novelty of the current \ARS.
We use the last remaining sufficient condition for well-balancedness from \eqref{eq:WB_consequence_on_intermediate_states} on $s$ to determine~$\delta E$.
Namely, we impose that, as soon as $(W_L, W_R)$ is an ISS,
the entropy of the intermediate states satisfies $s(W_L^*) = s(W_R^*) = s^*$.
Moreover, due to~\eqref{eq:condition_on_s_star}, we know that setting~$s^*$
according to \eqref{eq:condition_on_s_star} yields the entropy stability
under the prerequisite that an entropy-entropy flux pair exists.
Following the discussion from \cref{sec:equation_of_state},
the expression of the internal energy $e$ with respect to
density $\rho$ and entropy $s$ is available.
This allows us to define the intermediate internal energies as
\begin{equation*}
    e_L^* = e(\rho_L^*, s^*) \text{\qquad and \qquad} e_R^* = e(\rho_R^*, s^*).
\end{equation*}
{Recall that the internal energy, given by the EOS, is a positive function.
Thus, the intermediate internal energies
$e_L^*$ and $e_R^*$ are also positive,
and \cref{lem:positivity_inte} holds.}

Transforming back to the total energy, we obtain, using \eqref{eq:TotalEnergy},
\begin{equation}
    \label{eq:linear_system_for_delta_E}
    E_L^* = \rho_L^* e_L^* + \frac{{(q_L^*)^2}}{2\rho_L^*} = \widehat E - \delta E
     \text{\qquad and \qquad}
    E_R^* = \rho_R^* e_R^* + \frac{{(q_R^*)^2}}{2\rho_R^*} = \widehat E + \delta E.
\end{equation}
Since the momentum is constant for an ISS,
we replace {$(q_L^*)^2$ and $(q_R^*)^2$}
in the above expression with a parameter \smash{$\widetilde{q^2}$}.
Equation \eqref{eq:linear_system_for_delta_E} is nothing but
a linear system in \smash{$(\widetilde{q^2}, \delta E)$}, whose solution is
\begin{equation*}
    \widetilde{q^2} = \frac{2 \rho_L^* \rho_R^*}{\rho_L^* + \rho_R^*} \left( 2 \widehat E - \rho_L^* e_L^* - \rho_R^* e_R^*\right)
     \text{\qquad and \qquad}
    \delta E = \frac{1}{2} \left(\rho_R^* e_R^* - \rho_L^* e_L^*\right) + \frac{{\smash{\widetilde{q^2}}}}{2} \frac{\delta \rho}{2 \rho_L^* \rho_R^*}.
\end{equation*}
Since the derivation started by considering a constant entropy,
the above expressions are well-balanced by construction, for any EOS.
Moreover, note that the expressions of \smash{$\widetilde{q^2}$} and $\delta E$ derived here
are a generalization of the ideal gas case considered in \cite{Berthon2025},
and the above formulas reduce to the ideal gas case when setting $e = (\gamma - 1) \, p / \rho$.

\subsection{\texorpdfstring{Determination of $S^q$ and $S^E$}{Determination of Sq and SE}}
\label{sec:source_term_approximation}

With the expression for \(\delta W\) established,
we now focus on determining the remaining two unknowns, \(S^q\) and \(S^E\).
To do so, we make use of the expressions for \(\smash{\widehat W}\)
given in \eqref{eq:expression_w_hat}.
Additionally, we recall that if \((W_L, W_R)\) is an ISS,
then \(\widehat W = \xoverline W\),
where \(\xoverline X \coloneqq (X_L + X_R) / 2\)
represents the arithmetic mean of the two quantities \(X_L\) and \(X_R\).
Substituting the expression of \(W_\HLL\) into~\eqref{eq:expression_w_hat},
we find that \(S^q\) and \(S^E\) must satisfy the following conditions
when \((W_L, W_R)\) is an~ISS:
\begin{subequations}
    \label{eq:conditions_on_Sq_and_SE}
    \begin{align}
        \label{eq:conditions_on_Sq}
        S^q \dx & =
        \left(\frac{q_R^2}{\rho_R} + p_R\right)
        -
        \left(\frac{q_L^2}{\rho_L} + p_L\right), \\
        \label{eq:conditions_on_SE}
        S^E \dx & =
        \left(\frac{q_R}{\rho_R} (E_R + p_R) \right)
        -
        \left(\frac{q_L}{\rho_L} (E_L + p_L)\right).
    \end{align}
\end{subequations}
The goal of this section is to find consistent approximations $S^q$ and $S^E$ of the source term averages satisfying \eqref{eq:conditions_on_Sq_and_SE}.


In contrast to $S^q$, the expression of the energy source term $S^E$ does not depend on the EOS.
Thus, we directly use the one from \cite{Berthon2025}, which is given by
\begin{equation}
    \label{eq:expression_of_SE}
    S^E = - \frac{q_L + q_R} 2 \, \frac{\varphi_R - \varphi_L} \dx.
\end{equation}
This expression is consistent with the source term $S^{(3)}(W) = -q \pd_x \varphi$ and satisfies the well-balanced condition \eqref{eq:conditions_on_SE}.

Next, we derive an approximation $S^q$ of the momentum source term such that the ISS condition \eqref{eq:conditions_on_Sq} is satisfied.
This approximation also depends on the considered EOS.
We make the intuitive Ansatz
\begin{equation}\label{eq:Ansatz_Sq}
    S^q \Delta x = -\frac{2 \rho_L \rho_R}{\rho_L + \rho_R} \jump{\varphi} + \varepsilon
\end{equation}
where $\varepsilon$ is a corrective term that accounts for the fact that the first consistent term does not fulfil the ISS condition \eqref{eq:conditions_on_Sq}.
However, away from the equilibrium state $\varepsilon$ should be consistent with zero.
To determine the correction~$\varepsilon$, we assume that $(W_L,W_R)$ form an ISS.
As we have seen in the derivation of~$\delta W$, in that case we have the relation \eqref{eq:WB_h} connecting $\jump{\varphi}$ and the enthalpy jump $\jump{h}$, i.e.,
\begin{equation}
    \label{eq:jump_phi_for_ISS}
    \jump{\varphi} = - \jump{e + \frac{p}{\rho} + \frac{1}{2}\frac{q^2}{\rho^2}}
\end{equation}
holds.
Thus, the dependence on the EOS needs to be taken into account and the correction term~$\varepsilon$ will differ from the one obtained in \cite{Berthon2025} for the ideal gas law.
Equating the right-hand sides of \eqref{eq:Ansatz_Sq} and \eqref{eq:conditions_on_Sq},
and using \eqref{eq:jump_phi_for_ISS},
a lengthy but straightforward calculation yields
\begin{equation}\label{eq:Sq_eps}
    \varepsilon
    =
    - \frac{2 \rho_L \rho_R}{\rho_L + \rho_R}
    \biggl(
    e(\rho_{R}, \xoverline{s}) - e(\rho_{L}, \xoverline{s})
    +
    \frac{p(\rho_{L}, s_L) + p(\rho_{R}, s_R)} 2
    \biggl( \frac{1}{\rho_R} - \frac{1}{\rho_L} \biggr)
    \biggr).
\end{equation}
As defined above, $\xoverline{s} = (s_L + s_R)/2$ is the arithmetic mean of the left and right entropy states.
Since $\varepsilon$ is not consistent with zero within a shock wave,
we proceed as in \cite{Berthon2025} and multiply it
with the function $\psi$ from \cref{lem:helper_function_for_density} which ensures that the resulting discretization of $S^q$ is consistent with the source term $S^{(2)}(W)$, while preserving the well-balanced property.
This is summarized in the following result.

\begin{lemma}
    Let $W_L, W_R \in \Omega$ be two given states and $S^q$ be defined as
    \begin{equation}
        \label{eq:expression_of_Sq}
        S^q
        =
        - \frac{2 \rho_L \rho_R}{\rho_L + \rho_R} \frac{\jump{\varphi}}{\dx}
        + \frac{\varepsilon}{\dx} \, \psi\big( [\varphi], [h], 3\big),
    \end{equation}
    with $\varepsilon$ given by \eqref{eq:Sq_eps}.
    Then $S^q$ is consistent with $S^{(2)}(W) = -\rho \pd_x \varphi$
    and satisfies the well-balanced condition \eqref{eq:conditions_on_Sq}.
\end{lemma}
\begin{proof}
    The well-balanced property is satisfied by construction,
    since $\psi = 1$ as soon as $[\varphi] = -[h]$.
    It is left to show the consistency.
    In particular, note that the first term in \eqref{eq:Ansatz_Sq} is consistent with $-\rho \pd_x \varphi$ as it is a first-order discretization of the source term.
    Therefore, we need to show that
    \begin{equation}
        \label{eq:consistency_requirement_eps}
        \lim_{\dx \to 0}
        \frac{\varepsilon}{\dx} \, \psi\big( [\varphi], [h], 3\big) = 0.
    \end{equation}
    Note that $\varphi$ is a given smooth function, i.e. $\jump{\varphi} = \mathcal{O}(\Delta x)$.
    We need to distinguish two cases:
    (a) the solution is smooth, and (b) the solution is not smooth.

    \begin{enumerate}[nosep]

        \item[(a)] We consider a smooth solution,
              i.e., $\rho$ and $s$ are smooth functions of $x$.
              We define the left and right densities as
              \begin{equation*}
                  \rho_L = \rho(x), \quad \rho_R = \rho(x+\Delta x).
              \end{equation*}
              Analogously, we define the left and right entropies.
              Then, Taylor series expansions yield
              \begin{equation*}
                  \frac{e(\rho_{R}, \xoverline{s}) - e(\rho_{L}, \xoverline{s})}{\dx}
                  =
                  {\left(\frac{\partial e}{\partial \rho}\right)_{\!\!s\!}} (\rho(x),\xoverline{s}) \, \partial_x \rho(x) + \mathcal{O}(\dx),
              \end{equation*}
              {where the subscript $s$ denotes the derivative at constant entropy,}
              and
              \begin{equation*}
                  \frac 1 \dx \jump{\frac{1}{\rho}}
                  =
                  - \frac{\partial_x \rho(x)}{\rho(x)^2} + \mathcal{O}(\Delta x).
              \end{equation*}
              Using {the harmonic mean
                      \begin{equation*}
                          \frac{2 \rho_L \rho_R}{\rho_L + \rho_R} = \rho(x) + \mathcal{O}(\Delta x)
                      \end{equation*}
                      and the arithmetic mean
                      \begin{equation*}
                          \frac{p(\rho_{L}, s_L) + p(\rho_{R}, s_R)} 2
                          =
                          p(\rho(x),\xoverline{s})
                          + \mathcal{O}(\Delta x),
                      \end{equation*}}
              we obtain
              \begin{equation}
                  \label{eq:expression_of_epsilon_over_dx_before_Gibbs}
                  \frac{\varepsilon}{\Delta x} =
                  - \rho \, \partial_x \rho \left({\left(\frac{\partial e}{\partial \rho}\right)_{\!\!s\!}} - \frac{p}{\rho^2} \right) + \mathcal{O}(\Delta x).
              \end{equation}
              {Now, recall from Gibbs' relation \eqref{eq:DerivEntropy} that
              \begin{equation*}
                  \left(\frac{\partial e}{\partial \rho}\right)_{\!\!s\!} = \frac{p}{\rho^2},
              \end{equation*}
              from which it follows that ${\varepsilon} / {\Delta x} = \mathcal{O}(\Delta x)$ in \eqref{eq:expression_of_epsilon_over_dx_before_Gibbs}.}
              Lengthy computations also show that the second-order term in the Taylor expansion vanishes, and we obtain ${\varepsilon} / {\Delta x} = \mathcal{O}(\Delta x^2)$.

        \item[(b)] Let us now turn to the case where the solution is not smooth,
              i.e., $W_R = W_L + \mathcal{O}(1)$.
              In this case, we have $\varepsilon = \mathcal{O}(1)$,
              and \eqref{eq:psi_vi} ensures that
              \smash{$\psi\big( [\varphi], [h], 3\big) = \mathcal{O}([\varphi]^3)$}.
              Since $\varphi$ is smooth, we have $[\varphi] = \mathcal{O}(\dx)$,
              and so ${\varepsilon} / {\Delta x} = \mathcal{O}(\Delta x^2)$.

    \end{enumerate}
    Therefore, in both cases, ${\varepsilon} / {\Delta x} = \mathcal{O}(\Delta x^2)$.
    This means that the consistency requirement \eqref{eq:consistency_requirement_eps}
    is satisfied in both cases, which concludes the proof.
\end{proof}

\subsection{Summary of the numerical scheme and main properties}

We state the fully defined approximate Riemann solver
and summarize the main properties of the associated Godunov-type finite volume scheme.
The \ARS \eqref{eq:ARS} is given by
\begin{equation}
    \label{eq:ARS_recap}
    \widetilde{W}\left( \frac x t; W_L, W_R \right) =
    \begin{dcases}
        W_L                   & \text{if} \quad\hskip15mm x < -\lambda t,   \\
        \widehat W - \delta W & \text{if} \quad -\lambda t < x < 0,         \\
        \widehat W + \delta W & \text{if} \quad\hskip7mm 0 < x < \lambda t, \\
        W_R                   & \text{if} \quad\hskip14.5mm x > \lambda t,  \\
    \end{dcases}
\end{equation}
where the components of $\widehat W$ are given by
\begin{equation}
    \label{eq:expression_w_hat_recap}
    \left\{
    \begin{aligned}
        \widehat \rho & = \rho_\HLL,   \\
        \widehat q    & = q_\HLL
        + \frac {S^q \dx} {2 \lambda}, \\
        \widehat E    & = E_\HLL
        + \frac {S^E \dx} {2 \lambda}, \\
    \end{aligned}
    \right.
\end{equation}
with $W_\HLL$ given by \eqref{eq:def_w_hll}
and the approximate wave speed $\lambda$ satisfying \eqref{eq:def_lambda}.
In \eqref{eq:expression_w_hat_recap},
the source term approximations $S^q$ and $S^E${, which vanish when $\varphi_L = \varphi_R$,}
are respectively given by
\eqref{eq:expression_of_Sq} and \eqref{eq:expression_of_SE}, i.e.,
\begin{equation}
    \label{eq:expressions_of_Sq_and_SE_recap}
    \begin{aligned}
        S^q
            & =
        - \frac{2 \rho_L \rho_R}{\rho_L + \rho_R} \frac{\varphi_R - \varphi_L} \dx
        + \frac{\varepsilon}{\dx} \, \psi\big( \varphi_R - \varphi_L, h_R - h_L, 3\big), \\
        S^E & = - \frac{q_L + q_R} 2 \, \frac{\varphi_R - \varphi_L} \dx,
    \end{aligned}
\end{equation}
where $\psi$ is defined in \cref{lem:helper_function_for_density},
and where $\varepsilon$ satisfies, according to \eqref{eq:Sq_eps},
\begin{equation*}
    \label{eq:expression_of_epsilon_recap}
    \varepsilon
    =
    - \frac{2 \rho_L \rho_R}{\rho_L + \rho_R}
    \biggl(
    e(\rho_{R}, \xoverline{s}) - e(\rho_{L}, \xoverline{s})
    +
    \frac{p(\rho_{L}, s_L) + p(\rho_{R}, s_R)} 2
    \biggl( \frac{1}{\rho_R} - \frac{1}{\rho_L} \biggr)
    \biggr),
\end{equation*}
with $\xoverline s = (s_L + s_R) / 2$.
To define $\delta W$, we need the intermediate quantities
\begin{equation*}
    s^* = \frac{(\rho s)_\HLL}{\rho_\HLL}
    \text{, \qquad}
    e_L^* = e(\rho_L^*, s^*)
    \text{\qquad and \qquad}
    e_R^* = e(\rho_R^*, s^*).
\end{equation*}
The components of $\delta W$ are then constructed in \cref{sec:determination_of_delta_w},
and they satisfy
\begin{equation}
    \label{eq:expression_delta_w_recap}
    \left\{
    \begin{aligned}
        \delta \rho &
        = \frac {\rho_R - \rho_L} 2 \,
        \psi \big( \varphi_R - \varphi_L, h_R - h_L, 1 \big),         \\
        \delta q    & = 0 \vphantom{\dfrac 1 2},                      \\
        \delta E    & =
        \frac{1}{2} \left(\rho_R^* e_R^* - \rho_L^* e_L^*\right)
        + \frac 1 4 \frac{\rho_R^* - \rho_L^*} {\rho_L^* + \rho_R^*}
        \left( 2 \widehat E - \rho_L^* e_L^* - \rho_R^* e_R^*\right). \\
    \end{aligned}
    \right.
\end{equation}

The properties of the numerical scheme are summarized in the main result.

\begin{theorem}
    \label{theo:summary}
    Let the time step $\dt$ be given by \eqref{eq:CFL_condition}
    and assume that the initial {data $W_i^0$ satisfy}
    $W_i^0 \in \Omega$ for all $i \in \mathbb{Z}$.
    Then, the numerical scheme \eqref{eq:scheme_from_ARS}
    with the approximate Riemann solver \eqref{eq:ARS_recap},
    where $\smash{\widehat{W}}$ is given by \eqref{eq:expression_w_hat_recap}
    and $\delta W$ is given by \eqref{eq:expression_delta_w_recap},
    satisfies the following properties:
    \begin{enumerate}
        \item \underline{\smash{consistency}} with the Euler equations with gravity \eqref{eq:EulerG};
        \item \underline{\smash{positivity}} of the density {and internal energy:}
              for all $n \geq 0$,
              \begin{equation*}
                  \forall i \in \mathbb{Z}, \; W_i^{n} \in \Omega
                  \quad \implies \quad
                  \forall i \in \mathbb{Z}, \; W_i^{n+1} \in \Omega;
              \end{equation*}
        \item \underline{\smash{entropy stability}}:
              {for all mathematical entropies $\rho \eta(s)$,
              for all $i \in \mathbb{Z}$, and for all $n \geq 0$,}
              \begin{equation*}
                  \rho_i^{n+1}\eta(s_i^{n+1}) \leqslant \rho_i^n \eta(s_i^n)
                  - \frac{\dt}{\dx} \left((\rho \eta(s) u)_\iph^n -
                  (\rho \eta(s) u)_\imh^n\right);
              \end{equation*}
        \item \underline{\smash{well-balancedness}}:
              \begin{equation*}
                  \forall i \in \mathbb{Z}, \; (W_i^n, W_{i+1}^n) \text{ is an ISS}
                  \quad \implies \quad
                  \forall i \in \mathbb{Z}, \; W_i^{n+1} = W_i^n.
              \end{equation*}
    \end{enumerate}
\end{theorem}

\begin{proof}
    We prove the four properties in order,
    using the results derived in the previous sections.
    \begin{enumerate}
        \item According to \cite{HarLaxLee1983},
              the scheme is consistent as soon as the \ARS satisfies
              the integral consistency relation,
              see \cref{sec:Riemann_solver_consistency}.
              The \ARS has been constructed such that this relation holds,
              and thus it is consistent.
              Therefore, the numerical scheme \eqref{eq:scheme_from_ARS}
              is also consistent.
        \item The positivity of the density is guaranteed by
              that of the intermediate densities.
              Indeed, according to \cref{lem:positivity_density},
              the intermediate densities satisfy
              the positivity condition \eqref{eq:positivity_condition_on_rho}.
              Therefore, $\rho_L^\ast > 0$ and $\rho_R^\ast > 0$ in the \ARS,
              and thus $\rho_i^{n+1} > 0$ in the numerical scheme.
                  {Further, \cref{lem:positivity_inte} yields
                      the positivity of the internal energy, i.e., $e_i^{n+1} > 0$,
                      since the intermediate internal energies $e_L^*$ and $e_R^*$ are positive.}
        \item {The entropy stability of the scheme is
              ensured since the
              intermediate specific entropy satisfies
              \eqref{eq:condition_on_s_star}.
              Indeed, under this condition,
              the \ARS \eqref{eq:ARS_recap}
              fulfills the entropy consistency condition
              \eqref{eq:Entropyconsistency}, which,
              following \cite{Berthon2025},
              proves that the scheme is entropy-stable.}
        \item \cref{cor:WB_consequence} gives
              a sufficient condition for well-balancedness,
              which leads to the conditions
              \eqref{eq:WB_consequence_on_hat_and_delta_w}
              on $\smash{\widehat{W}}$ and $\delta W$
              used in \eqref{eq:ARS_recap} to define the \ARS.
              \cref{sec:determination_of_delta_w,sec:source_term_approximation}
              constructs $\smash{\widehat{W}}$ and $\delta W$
              such that the conditions \eqref{eq:WB_consequence_on_hat_and_delta_w}
              are satisfied,
              and consequently the numerical scheme is well-balanced.
    \end{enumerate}
    All four properties have been proven, and the proof is thus concluded.
\end{proof}

{
\subsection{Second-order extension}
Instead of considering the numerical scheme \eqref{eq:ARS} in terms of Riemann solvers, it is more convenient, for the second-order extension, to start from the scheme in standard form
\begin{equation}
    \label{eq:schemeStandard}
    W_i^{n+1} = W_i^n
    -
    \frac{\dt}{\dx}\left(
    \mathcal{F}(W_\iphm^{n},W_\iphp^{n})
    - \mathcal{F}(W_\imhm^{n,}, W_\imhp^{n})
    \right) + \Delta t \overline{\mathcal{S}}_i^n.
\end{equation}
Therein, $\mathcal{F}$ denotes the numerical flux, which is evaluated at the cell interface values $W_\imhpminline^n$ and $W_\iphpminline^n$ at interfaces $x_\imhinline$ and $x_\iphinline$, respectively.
Following the derivation given in e.g. \cite{HarLaxLee1983}, the numerical flux for two input values $W_L$ and $W_R$ takes the form
\begin{equation*}
    \mathcal{F}(W_L,W_R) = \frac{F(W_L)+F(W_R)}{2} - \frac{\lambda}{2}\left(W_L^\ast - W_L\right) + \frac{\lambda}{2}\left(W_R^\ast - W_R\right).
\end{equation*}
Following \cite{Desveaux2022}, the numerical source term is given by
\begin{equation}\label{eq:sourceTrap}
    \overline{\mathcal{S}}_i
    =
    \frac{1}{2}\left(
    \mathcal{S}(W_\imhm^{n,}, W_\imhp^{n})
    + 2 \mathcal{S}(W_\imhp^n, W_\iphm^n)
    + \mathcal{S}(W_\iphm^{n},W_\iphp^{n})
    \right)
\end{equation}
where
\begin{equation*}
    \mathcal{S}(W_L,W_R) = \begin{pmatrix}
        0 \\ S^q \\ S^E
    \end{pmatrix},
\end{equation*}
as defined in \eqref{eq:expressions_of_Sq_and_SE_recap}.}

{For the first-order scheme derived in the previous sections, the interface values are given by $\smash{W_\imhminline^n = W_{i-1}^n}$, $\smash{W_\imhpinline^n = W_i^n}$, $\smash{W_\iphminline^n = W_i^n}$ and $\smash{W_\iphpinline^n = W_{i+1}^n}$, which corresponds to a constant reconstruction and a trapezoidal rule for the gravitational source terms.
To obtain a second-order reconstruction while ensuring the well-balanced property, we consider piecewise linear functions in so-called equilibrium variables $W^{\eqvar} = (q, s, H)^\intercal$.
Thus, we consider on the cell $C_i$ a linear function
\begin{equation} \label{eq:lin_recon}
    r_i^\eqvar(x) = W_i^\eqvar + \sigma_i (x - x_i),
\end{equation}
where $\sigma_i = (\sigma_i^q, \sigma_i^s, \sigma_i^H)^\intercal$ denotes the slope for each variable.
The slopes are computed from the neighbouring cells by
\begin{equation}\label{eq:slope}
    \sigma_i = \texttt{limiter}(W_{i-1}^\eqvar, W_i^\eqvar, W_{i+1}^\eqvar ),
\end{equation}
where the \texttt{limiter} is a consistent limiter function such as \texttt{minmod}, see e.g.~\cite{Lee1979}.
The limiter is called consistent if $W_{i-1}^\eqvar = W_i^\eqvar = W_{i+1}^\eqvar$ implies that $\sigma_i = 0$.
The interface values on cell $C_i$ are then computed by
\begin{equation}\label{eq:recon_val}
    W_\imhp^{\eqvar} = W_i^\eqvar - \sigma_i \frac{\Delta x}{2}
    \text{\qquad and \qquad}
    W_\iphm^{\eqvar} = W_i^\eqvar + \sigma_i \frac{\Delta x}{2}.
\end{equation}
The state variables at the interface $W_{i\pm1/2,\mp}^n$ are then recovered from the equilibrium ones by inverting the mapping $W \mapsto W^\eqvar$
(should this mapping not be invertible, one can instead use the classical reconstruction in primitive variables).
The given gravitational potential is reconstructed at the cell interfaces.
We are now able to state the first result of this section, regarding the reconstructed values.
\begin{lemma}\label{lem:slopes}
    Let the pairs $(W_{i-1}^n, W_i^n)$ and $(W_i^n,W_{i+1}^n)$ be ISSs. Then after applying the reconstruction \eqref{eq:lin_recon} -- \eqref{eq:recon_val} with a consistent limiter in \eqref{eq:slope}, we obtain $W_\imhpinline^n = W_i^n$ and $W_\iphminline^n = W_i^n$.
\end{lemma}
\begin{proof}
    By \cref{def:ISS}, the ISS assumption implies that $q^n_{i-1} = q^n_i = q^n_{i+1}$, $s^n_{i-1} = s^n_i = s^n_{i+1}$ and $H^n_{i-1} = H^n_i = H^n_{i+1}$. Thus, $\sigma_i = 0$ since the limiter is consistent. This implies $W_\imhpinline^{\eqvar} = W_\iphminline^{\eqvar} = W_i^\eqvar$, and thus $W_\imhpinline^n = W_i^n$ and $W_\iphminline^n = W_i^n$.
\end{proof}
}
Note that by transforming the reconstructed equilibrium variables into state variables, the positivity of density and internal energy might be violated. Hence, \cref{lem:positivity_density,lem:positivity_inte} would no longer apply, since they require $W_L, W_R \in \Omega$ and thus $W_\iphpminline^n \in \Omega$.
In \cite{Berthon2005, Thomann2019} explicit formulas for a slope limiting have been formulated for the reconstruction in primitive variables.
However, since here the transform from equilibrium to state variables do not have an analytical expression, we cannot give an explicit expression.
Therefore, we propose to reduce the slope adaptively by successively halving the slopes of the state variables until $W_\iphpminline^n \in \Omega$.
Note that this procedure does not violate the result from \cref{lem:slopes}.

We have shown how to obtain second-order accuracy in space.
To obtain second-order accuracy also in time, we first rewrite the scheme \eqref{eq:schemeStandard}, which is first-order in time and second-order in space, under the compact form
\begin{equation*}
    W_i^{n+1} = W_i^n - \Delta t \mathcal{H}_i^n,
\end{equation*}
where $\mathcal{H}_i^n$ denotes the spatial operator
\begin{equation*}
    \mathcal{H}_i^n = \frac{1}{\dx}\left(
    \mathcal{F}(W_\iphm^{n},W_\iphp^{n})
    - \mathcal{F}(W_\imhm^{n,}, W_\imhp^{n})
    \right) + \overline{\mathcal{S}}_i^n.
\end{equation*}
Then, we apply the two-stage time integration derived in \cite{Berthon2005}, which has the advantage that at each stage the CFL criterion is met and the final update is a convex combination of true first-order schemes in time.
Denoting by $W_i^{(1)}$ and $W_i^{(2)}$ the stage values, we have, for all $i \in \mathbb{Z}$,
\begin{align}\label{eq:scheme2}
    \begin{split}
        W_i^{(1)} & = W_i^n - \Delta t_1  \mathcal{H}_i^n,         \vphantom{\dfrac 1 2}\\
        W_i^{(2)} & = W_i^{(1)} - \Delta t_2  \mathcal{H}_i^{(1)}, \vphantom{\dfrac 1 2}\\
        W_i^{n+1} & = (1-\kappa) W_i^n + \kappa W_i^{(2)},
        \text{\qquad with \quad} \kappa = \frac{\Delta t_1 \Delta t_2}{(\Delta t_1 + \Delta t_2)^2}.
    \end{split}
\end{align}
The total time increment is then given by
\begin{equation*}
    \Delta t = \frac{2 \Delta t_1 \Delta t_2}{\Delta t_1 + \Delta t_2}.
\end{equation*}
Note that, if $\Delta t_1 = \Delta t_2$, the well-known Strong Stability Preserving Runge-Kutta 2 (SSP-RK2) method is recovered.
We conclude the section with a result for the second-order scheme.
\begin{theorem}
    \label{theo:order2}
    Let the time step $\dt$ be given by
    \begin{equation*}\label{eq:CFL2}
        \Delta t  \leq \frac{1}{4} \frac{\Delta x}{\max_i\big(\lambda(W_\imhp^n, W_\iphm^n)\big)}.
    \end{equation*}
    and assume that the initial {data $W_i^0$ satisfy}
    $W_i^0 \in \Omega$ for all $i \in \mathbb{Z}$.
    Then, the second-order numerical scheme \eqref{eq:scheme2} with the reconstruction \eqref{eq:lin_recon} and adaptive slope limiting on the state variables satisfies the following properties:
    \begin{enumerate}
        \item \underline{\smash{consistency}} with the Euler equations with gravity \eqref{eq:EulerG};
        \item \underline{\smash{positivity}} of the density {and internal energy:}
              for all $n \geq 0$,
              \begin{equation*}
                  \forall i \in \mathbb{Z}, \; W_i^{n} \in \Omega
                  \quad \implies \quad
                  \forall i \in \mathbb{Z}, \; W_i^{n+1} \in \Omega;
              \end{equation*}
        \item \underline{\smash{entropy stability}}:
              {for all mathematical entropies $\rho \eta(s)$,
              for all $i \in \mathbb{Z}$, and for all $n \geq 0$,}
              \begin{equation*}
                  \rho_i^{n+1}\eta(s_i^{n+1}) \leqslant \rho_i^n \eta(s_i^n)
                  - \frac{\dt}{\dx} \left((\rho \eta(s) u)_\iph^n -
                  (\rho \eta(s) u)_\imh^n\right);
              \end{equation*}
        \item \underline{\smash{well-balancedness}}:
              \begin{equation*}
                  \forall i \in \mathbb{Z}, \; (W_i^n, W_{i+1}^n) \text{ is an ISS}
                  \quad \implies \quad
                  \forall i \in \mathbb{Z}, \; W_i^{n+1} = W_i^n.
              \end{equation*}
    \end{enumerate}
\end{theorem}
\begin{proof}
    The consistency follows from the consistency of the first-order scheme.
    \newline
    The second-order scheme  \eqref{eq:scheme2} with the source term \eqref{eq:sourceTrap} can be written as a convex combination of first-order schemes \eqref{eq:schemeStandard} for $W_\imhpinline^n$ and $W_\iphminline^n$.
    Due to the slope limiting on the state variables, the positivity of density and internal energy for the first-order in time and second-order  in space scheme  \eqref{eq:schemeStandard} follows from \cite[Theorem 2.1]{Berthon2005}. The entropy stability follows from \cite[Theorem 2.2]{Berthon2005} since by Theorem \ref{theo:summary} the first-order scheme fulfils properties 2 and 3.
    The fully second-order scheme \eqref{eq:scheme2} is a convex combination of first-order steps \eqref{eq:schemeStandard} which ensures $W_i^{n+1} \in \Omega$ for all $i \in \mathbb{Z}$ and the entropy stability.
    \newline
    Moreover, the well-balanced property is ensured since by \cref{lem:slopes} the space discretization reduces to the first-order one, and thus it follows that \smash{$W_i^{(1)} = W_i^n$} and \smash{$W_i^{(2)} = W_i^n$} as soon as a steady solution is considered. Therefore, by definition of the second-order time integrator, we obtain $W_i^{n+1} = W_i^n$.
\end{proof}

\section{Numerical results}
\label{sec:NumRes}

In this section, we perform numerical test cases to
validate the theoretical properties of the
first-order well-balanced scheme described
in \cref{theo:summary}.
{We denote by \FWBscheme the presented {first-order} fully well-balanced scheme, {by \FWBschemeTwo the second-order extension}
and by \HLLscheme the reference, non-well-balanced HLL scheme
from~\cite{HarLaxLee1983} with a centered discretization of the source terms.}

In all test cases, we apply six different EOSs to assess the performance of the numerical schemes.
Four applied EOS are obtained from the cubic EOS
\eqref{eq:EOS_cubic},
namely the ideal gas law (denoted by \ideal),
the van-der-Waals gas (denoted by \vdW),
the Redlich-Kwong EOS (denoted by \RedKwo),
and the Peng-Robinson EOS (denoted by \PenRob).
The free parameters of these four cubic EOSs
are given in \cref{tab:EOS_parameters}.
The remaining two EOSs are tabulated and are obtained by using the \texttt{CoolProp}
library~\cite{CoolProp}.
The first tabulated EOS, denoted by \water, models the properties of water,
while the second one, denoted by \methane, includes the physical properties of methane.

\begin{table}[!t]
    \centering
    \begin{tabular}{ccccc}
        \toprule
                 & \ideal      & \vdW                & \RedKwo & \PenRob \\
        \cmidrule(lr){1-5}
        $R \, [\unit{\J\per\K}]$
                 & $0.4$       & $0.4$               & $0.4$   & $0.4$   \\
        $c_v^0 \, [\unit{\J\per\K}]$
                 & $1$         & $1$                 & $1$     & $1$     \\
        $s_0 \, [\unit{\J\per\K}]$
                 & $\log(0.4)$ & $\log(0.4)$         & $0$     & $0$     \\
        $a_0 \, [^\text{(a)}]$
                 & ---         & $15.67^\text{(b)}$  & $15$    & $15$    \\
        $b \, [\unit{\m\cubed\per\kg}]$
                 & ---         & $0.1273^\text{(b)}$ & $0.05$  & $0.05$  \\
        $T_0 \, [\unit{\K}]$
                 & ---         & ---                 & ---     & $0.3$   \\
        $\kappa$ & ---
                 & ---         & ---                 & $0.5$             \\
        \bottomrule
    \end{tabular}
    \caption{Parameters for the ideal gas, van-der-Waals,
        Redlich-Kwong and Peng-Robinson EOS.
        A dash indicates that the parameter is not used.
        Notes:
        (a) the unit of $a_0$ is $\unit{\m^6\K^{1/2}\per\kg\squared}$ for \RedKwo
        and $\unit{\m^6\per\kg\squared}$ for the other three EOSs;
        (b) the values of $a_0$ and $b$ for the van-der-Waals EOS
        are modified in \cref{sec:homogeneous_RP}.
    }
    \label{tab:EOS_parameters}
\end{table}

Further, if not otherwise specified,
we apply in all test cases the quadratic gravitational potential
\begin{equation}
    \label{eq:default_phi_expression}
    \varphi(x) = \frac {\varphi_0} 2 (x-x_0)^2.
\end{equation}
The scaling factor $\varphi_0$ ensures that the flow is within the same Mach and Froude number regime.
Therefore, $\varphi_0$ depends on the typical magnitude of the equilibrium pressure and density.
We set $\varphi_0 = \qty{2e4}{\per\square\s}$ for both tabulated EOSs due to the physical properties of water and methane at normal pressure,  i.e. Mach numbers around $10^{-2}$,
and $\varphi_0 = \qty{1}{\per\square\s}$ for the four analytical EOSs which corresponds to Mach regimes around $1$.
For more details on well-balanced methods for low Mach and Froude number regimes and formal Mach and Froude number limits for the Euler equations with gravity, see e.g. \cite{Thomann2020}.

The time step is given by the equality cases of~\eqref{eq:CFL_condition}, {for \FWBscheme and \eqref{eq:CFL2} for \FWBschemeTwo.}
The parameter $\Lambda$
in the eigenvalue estimate~\eqref{eq:def_lambda} is
also used when computing the time step.
Unless otherwise mentioned, it is set to $1$.
    {Except for the accuracy test case where a linear reconstruction is used, the standard minmod function is used in \eqref{eq:slope} to compute the slopes. }
    {For none of the test cases, the slope limiting was triggered.}

\subsection{Accuracy of the numerical schemes}
\label{sec:accuracy}

\begin{table}[!t]
    \centering
    \begin{tabular}{ccccccc}
        \toprule
         & \ideal & \vdW   & \RedKwo & \PenRob       & \water         & \methane       \\
        \cmidrule(lr){1-7}
        $\rho_0 \, [\unit{\kg\per\m\cubed}]$
         & $2$    & $2$    & $2$     & $15$          & $\num{997.05}$ & $\num{422.8}$  \\
        $u_0 \, [\unit{\m\per\second}]$
         & $0.25$ & $0.25$ & $0.25$  & $0.25$        & $\num{1250}$   & $\num{1000}$   \\
        $p_0 \, [\unit{\Pa}]$
         & $5$    & $20$   & $5$     & $15$          & $\num{101800}$ & $\num{101800}$ \\
        $A$
         & $0.25$ & $0.1$  & $0.25$  & $\num{0.025}$ & $10^{-4}$      & $10^{-3}$      \\
        \bottomrule
    \end{tabular}
    \caption{Parameters of the exact solution \eqref{eq:exact_sol} from \cref{sec:accuracy}, for each of the six EOSs under consideration.}
    \label{tab:exact_solution_parameters}
\end{table}

To verify the experimental order of convergence (EOC) of the numerical scheme,
we study an exact solution of the Euler equations with gravity~\eqref{eq:EulerG}
taken from~\cite{KliPupSem2019},
which is a variation of the test introduced in \cite{XinShu2012}.
The analytical solution, with constant velocity $u_0$, is given by
\begin{equation}
    \label{eq:exact_sol}
    \begin{dcases}
        \rho(x,t)
        =
        \vphantom{\dfrac12}\rho_0 \left(1 + A \sin\big(k \pi (x - u_0 t)\big)\right), \\
        u(x,t)
        =
        \vphantom{\dfrac12}u_0,                                                       \\
        p(x,t)
        =
        p_0 - \rho_0 \left((x - u_0 t) - \frac{A}{k \pi} \cos\big(k \pi (x - u_0 t)\big)\right).
    \end{dcases}
\end{equation}
We take $k = 4$, which yields a highly oscillatory solution.
The other parameters are given in
\cref{tab:exact_solution_parameters}
with respect to each used EOS.
Indeed, since all EOSs have different domains of validity,
the parameters differ for each EOS.
This test case is designed for a linear potential $\varphi(x) = x$.
The simulation is carried out on the computational domain $\Omega = [0,1]$,
on seven grids with $N = 16 \cdot 2^{j}$ cells, where $j \in \{ 0, \dots, 6 \}$ with exact boundary conditions.
The final time is $t_\text{f} = 10^{-5}\,\unit{\s}$ for both tabulated EOSs,
and $t_\text{f} = \qty{0.5}{\s}$ for the other four EOSs.
    {This discrepancy between final times, which will be repeated in the following test cases, is due to the different scales in wave speeds between the tabulated EOSs and the other ones. Namely, the tabulated EOSs have a much higher speed of sound than the analytical ones, which leads to a faster propagation of waves.}

\begin{figure}[!t]
    \centering
    \includegraphics[scale=1.25]{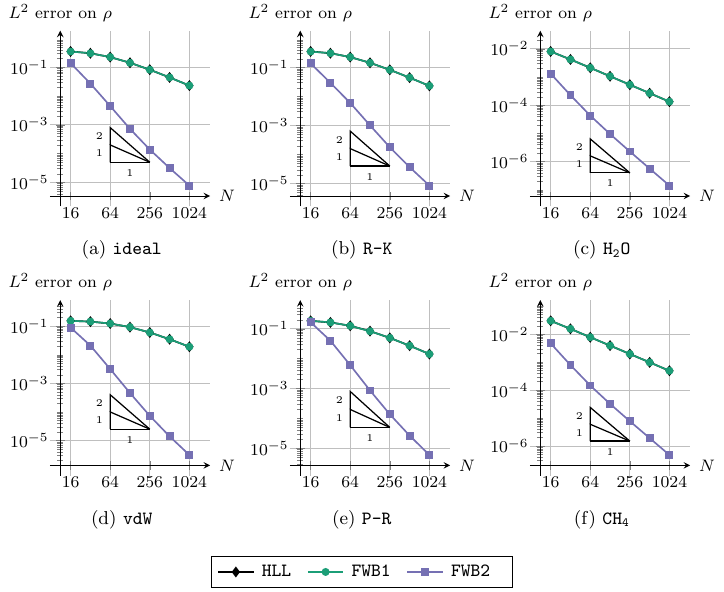}
    \caption{Experimental order of convergence and $L^2$ errors between the exact solution~\eqref{eq:exact_sol} and its approximation by the \HLLscheme, \FWBscheme and \FWBschemeTwo schemes, for each of the six EOSs under consideration. }
    \label{fig:EOC_phi1}
\end{figure}

In \cref{fig:EOC_phi1}, the $L^2$ errors are depicted
with respect to the number of discretization cells.
    {We recover the expected EOC {of one} for the \FWBscheme scheme {and two for the \FWBschemeTwo} independently of the considered EOS.
        Moreover, we note that, as expected, the results {for the first-order \FWBscheme scheme} are very close to the reference \HLLscheme approximation on this unsteady case.}

\begin{figure}[!t]
    \centering
    \includegraphics[scale=1.25]{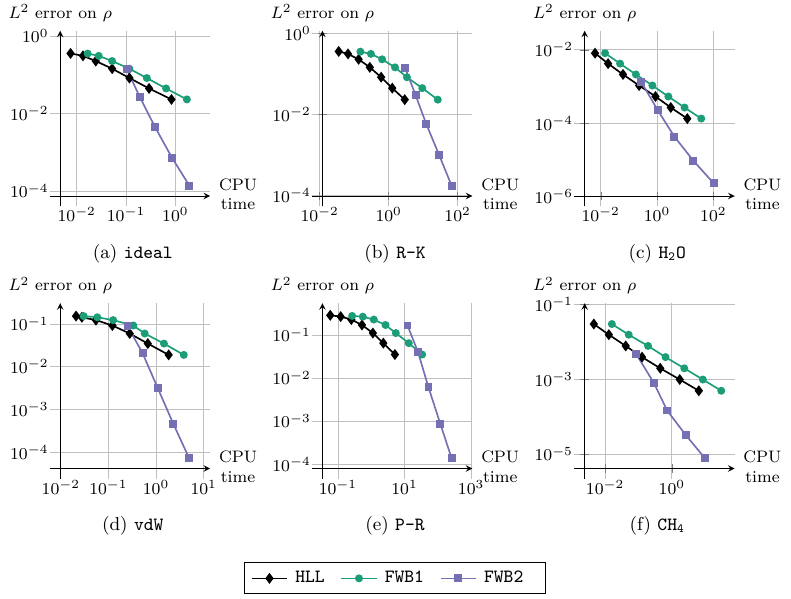}
    \caption{{Efficiency curves (error w.r.t. CPU time) for the three schemes, on the exact solution from \cref{sec:accuracy}, and for each of the six EOSs under consideration. Meshes with $({2^\ell})_{\ell \in \{4, \dots, 10\}}$ cells are used for the \HLLscheme and \FWBscheme schemes, and $({2^\ell})_{\ell \in \{4, \dots, 8\}}$ cells are considered for the \FWBschemeTwo scheme (more cells would yield very low errors, making it hard to plot the difference between the other two schemes).}}
    \label{fig:efficiency_exact_solution}
\end{figure}

    {We also compare, in \cref{fig:efficiency_exact_solution}, the efficiency of the three schemes on the approximation of this exact, unsteady solution. As expected, we observe that the \HLLscheme scheme is slightly more efficient than the \FWBscheme scheme, since the computational overhead caused by the well-balanced correction does not improve the error in this unsteady case. We also observe that this overhead is larger for EOSs such as \RedKwo and \PenRob where computing $e(\rho, s)$ is more expensive. However, the \FWBschemeTwo scheme is more efficient than the other two schemes, since it is second-order accurate in space and time.}

\subsection{Assessment of the well-balanced property}
\label{sec:numerical_moving}

The main motivation behind the construction of well-balanced schemes in general lies
in the resolution of small perturbations around equilibria using coarse meshes.
Applying a non-well-balanced scheme on those computational conditions easily leads to huge background errors stemming from the low resolution of the equilibrium state which makes the numerical solution unusable.
Thus, capturing in particular delicate perturbations would require a substantial grid refinement to reduce the truncation error and make the perturbations visible.

The next three sets of test cases therefore concern the performance of the novel fully well-balanced solvers on such flows on coarse grids.
The first one validates the exact preservation of a moving steady solution,
whereas the subsequent two are dedicated to studying perturbations of such equilibria.

The equilibrium solution is computed from the triplet $(q_0, s_0, H_0)$
given, with respect to each EOS, in \cref{tab:MovingEQ_Triplet}.
In all cases, the computational domain is $\Omega = [0,1]$ and the gravitational potential is given by \eqref{eq:default_phi_expression}.
The final times, for each experiment,
are reported in \cref{tab:MovingEQ_final_times} which differ for each EOS due to different flow regimes inflicted by initial condition and EOS.

\begin{table}[!bt]
    \centering
    \renewcommand{\arraystretch}{1.25}
    \begin{tabular}{rcccc}
        \toprule
        EOS      &
                 & $q_0 \left[\frac{\unit{\kg}}{\unit{\square\m\second}}\right]$
                 & $s_0 \left[\frac{\unit{\joule}}{\unit{\kelvin}}\right]$
                 & $H_0 \left[\frac{\smash{\unit{\square\m}}}{\unit{\square\second}}\right]$
        \\
        \cmidrule(lr){1-5}
        \ideal   &
                 & \num{1}
                 & \num{1}
                 & \num{5}
        \\
        \vdW     &
                 & \num{2.5}
                 & \num{-3.0}
                 & \num{55.0}
        \\
        \RedKwo  &
                 & \num{1}
                 & \num{-2.5}
                 & \num{12.5}
        \\
        \PenRob  &
                 & \num{5}
                 & \num{-2}
                 & \num{20}
        \\
        \water   &
                 & $10^5$
                 & \num{6.7e2}
                 & \num{2.1e5}
        \\
        \methane &
                 & $10^5$
                 & \num{-80}
                 & \num{2.5e4}
        \\
        \bottomrule
    \end{tabular}
    \caption{Values of the constant discharge $q_0$, entropy $s_0$ and specific total enthalpy $H_0$ for the moving equilibrium test cases from \cref{sec:numerical_moving}.}
    \label{tab:MovingEQ_Triplet}
\end{table}

\begin{table}[!bt]
    \centering
    \renewcommand{\arraystretch}{1.25}
    \begin{tabular}{rcccc}
        \toprule
        EOS      &
                 & $t_\text{f} \left[\unit{s}\right]$, \cref{sec:preservation_numerical}
                 & $t_\text{f} \left[\unit{s}\right]$, \cref{sec:perturbation_bump}
                 & $t_\text{f} \left[\unit{s}\right]$, \cref{sec:perturbation_sin}
        \\
        \cmidrule(lr){1-5}
        \ideal   &
                 & \num{1.5}
                 & \num{0.05}
                 & \num{0.72}
        \\
        \vdW     &
                 & \num{1.5}
                 & \num{0.02}
                 & \num{0.25}
        \\
        \RedKwo  &
                 & \num{1.5}
                 & \num{0.05}
                 & \num{0.72}
        \\
        \PenRob  &
                 & \num{0.1}
                 & $10^{-3}$
                 & \num{1.34}
        \\
        \water   &
                 & $10^{-4}$
                 & $10^{-4}$
                 & \num{5.6e-4}
        \\
        \methane &
                 & \num{2e-3}
                 & $10^{-4}$
                 & \num{7.2e-4}
        \\
        \bottomrule
    \end{tabular}
    \caption{Values of the final time $t_\text{f}$ for the moving equilibrium test cases from \cref{sec:numerical_moving}.}
    \label{tab:MovingEQ_final_times}
\end{table}

\subsubsection{Preservation of a moving steady solution}
\label{sec:preservation_numerical}

To numerically verify the well-balanced property of the schemes,
we first compute the time evolution of the equilibrium solution with nonzero velocity.
We take $50$ discretization cells,
and apply exact boundary conditions.
We expect the solution to remain exact,
up to machine precision,
when using our fully well-balanced scheme.

The $L^2$ errors at time $t = t_\text{f}$ are given in
\cref{tab:WB_moving_linear}.
All well-balanced schemes are able to preserve the moving equilibrium up to machine precision for all considered EOS,
whereas the classical HLL scheme \cite{HarLaxLee1983} yields quite large errors.
This verifies and illustrates the ability of the new well-balanced solvers
to capture moving equilibria up to machine precision for general EOS.

\begin{table}[tb!]
    \centering
    \begin{tabular}{ccccc}
        \toprule
                       &                &
        \HLLscheme     & \FWBscheme     & {\FWBschemeTwo}   \\
        \cmidrule(lr){1-5}
        \multirow{3}{*}{\ideal}
                       & $\rho$         &
        \num{3.94E-03} & \num{2.30E-15} & {\num{4.98e-17}}  \\
                       & $q$            &
        \num{3.38E-05} & \num{1.03E-16} & {\num{1.88e-14}}  \\
                       & $E$            &
        \num{3.38E-05} & \num{1.03E-16} & {\num{1.63e-16}}  \\
        \cmidrule(lr){1-5}
        \multirow{3}{*}{\vdW}
                       & $\rho$         &
        \num{3.97E-05} & \num{7.02E-16} & {\num{1.08e-15}}  \\
                       & $q$            &
        \num{6.29E-07} & \num{4.40E-17} & {\num{1.85e-15}}  \\
                       & $E$            &
        \num{2.18E-06} & \num{1.74E-16} & {\num{1.26e-15}}  \\
        \cmidrule(lr){1-5}
        \multirow{3}{*}{\RedKwo}
                       & $\rho$         &
        \num{2.15E-04} & \num{2.19E-15} & {\num{4.88e-15}}  \\
                       & $q$            &
        \num{8.76E-06} & \num{9.79E-17} & {\num{3.78e-15}}  \\
                       & $E$            &
        \num{2.87E-05} & \num{3.44E-16} & {\num{6.14e-15}}  \\
        \cmidrule(lr){1-5}
        \multirow{3}{*}{\PenRob}
                       & $\rho$         &
        \num{9.83E-05} & \num{1.95E-14} & {\num{1.35e-14} } \\
                       & $q$            &
        \num{4.40E-05} & \num{3.83E-15} & {\num{5.73e-15}}  \\
                       & $E$            &
        \num{2.68E-05} & \num{1.14E-14} & {\num{5.71e-15}}  \\
        \cmidrule(lr){1-5}
        \multirow{3}{*}{\water}
                       & $\rho$         &
        \num{1.13E-05} & \num{5.52E-14} & {\num{1.81e-14}}  \\
                       & $q$            &
        \num{8.07E-07} & \num{4.67E-14} & {\num{1.17e-13}}  \\
                       & $E$            &
        \num{1.16E-05} & \num{1.09E-13} & {\num{6.78e-14}}  \\
        \cmidrule(lr){1-5}
        \multirow{3}{*}{\methane}
                       & $\rho$         &
        \num{5.40E-06} & \num{1.46E-15} & {\num{2.00e-15}}  \\
                       & $q$            &
        \num{5.11E-05} & \num{1.70E-13} & {\num{3.74e-15}}  \\
                       & $E$            &
        \num{1.13E-04} & \num{4.68E-14} & {\num{1.27e-13}}  \\
        \bottomrule
    \end{tabular}
    \caption{Well-balanced test case from \cref{sec:numerical_moving}: $L^2$ errors on the density, the momentum and the total energy, reported for each of the six EOSs.}
    \label{tab:WB_moving_linear}
\end{table}

\subsubsection{Gaussian perturbation of a moving steady solution}
\label{sec:perturbation_bump}

Next, we consider a Gaussian perturbation of a moving equilibrium.
This is a standard test case whose purpose is to verify that underlying equilibrium does not inflict spurious errors in the perturbation.
The computational domain is discretized using $50$ cells, with inhomogeneous Dirichlet boundary conditions corresponding to the exact, unperturbed steady solution.
From the triplet $(q_0,s_0,H_0)$ reported in \cref{tab:MovingEQ_Triplet},
we compute the equilibrium state $(\rho_\text{eq}, u_\text{eq}, p_\text{eq})$.
Then, the initial density and velocity are set to $\rho(x,0) = \rho_\text{eq}(x)$ and $u(x,0) = u_\text{eq}$,
while a small perturbation is added to the pressure, i.e.,
\begin{equation*}
    p(x,0) = p_\text{eq}(x) \left(1 + \nu \exp\left(-100 \left( x - \frac{1}{2}\right)^2\right)\right).
\end{equation*}
The parameter $\nu$ denotes the amplitude of the initial perturbation,
and is set to $\nu = 10^{-4}$.
Different choices of $\nu$ yield similar numerical results and are thus omitted.
    {In this section and the following one, reference solutions are provided by running the proposed scheme with a fine mesh (made of \num{7500} cells). The perturbations are too small to be captured by traditional, non-well-balanced numerical schemes, thus a reference with the non-well-balanced \HLL~ scheme is omitted.}

\begin{figure}[tb!]
    \centering
    \includegraphics[scale=1.25]{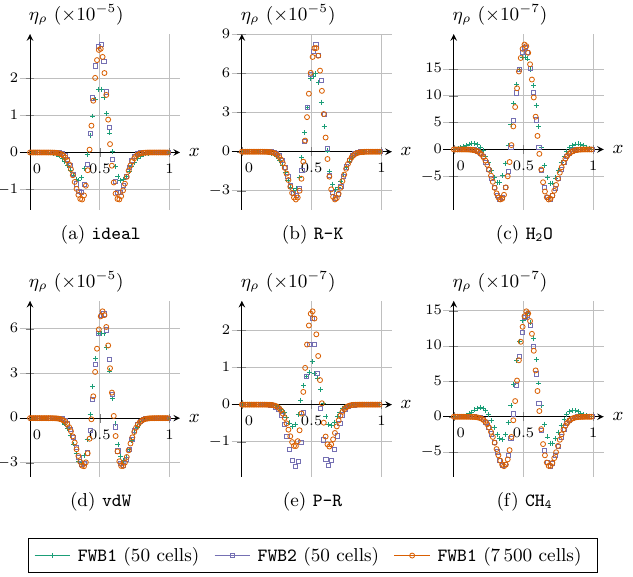}
    \caption{Perturbation of equilibrium state, described in \cref{sec:perturbation_bump}: Relative density difference $\reldiff_\rho = (\rho_\text{eq}(x) - \rho)/\rho_\text{eq}(x)$ with respect to the equilibrium density $\rho_\text{eq}(x)$, where $\rho$ is obtained by the \FWBscheme scheme with a mesh made of $50$ and \num{7500} discretization cells {and \FWBschemeTwo scheme on $50$ cells}, using the six EOSs under consideration.}
    \label{fig:moving_perturbed}
\end{figure}

In \cref{fig:moving_perturbed} we plot the relative density perturbation
$\reldiff_\rho = (\rho_\text{eq}(x) - \rho)/\rho_\text{eq}$,
scaled with respect to the background equilibrium $\rho_\text{eq}(x)$,
at the final time $t = t_\text{f}$.
The considered times guarantee that the waves triggered by the perturbation
are still contained inside the computational domain.
The perturbations for all considered EOS are well-captured
and no spurious errors are introduced from the background equilibrium.
    {Using the second-order well-balanced scheme \FWBschemeTwo on the coarse grid shows a great improvement, and results are comparable with the solution of the first-order scheme \FWBscheme on a very fine mesh.}

\begin{figure}[!t]
    \centering
    \includegraphics[scale=1.25]{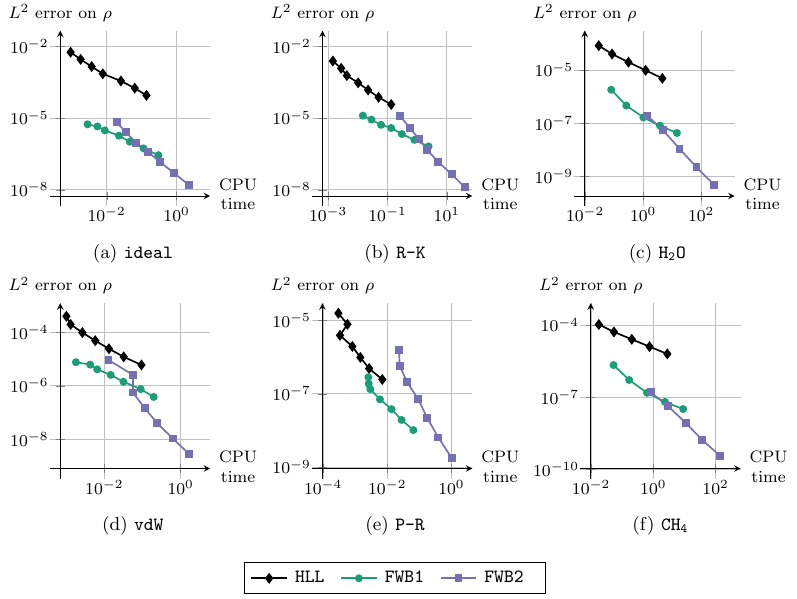}
    \caption{{Efficiency curves (error w.r.t. CPU time) for the three schemes, on the perturbed steady solution from \cref{sec:perturbation_bump}, and for each of the six EOSs under consideration. Meshes with $({2^\ell})_{\ell \in \{4, \dots, 10\}}$ cells are used for the \ideal, \vdW, \RedKwo and \PenRob EOSs, and $({2^\ell})_{\ell \in \{4, \dots, 8\}}$ cells are considered for the \water and \methane ones.}}
    \label{fig:efficiency_steady_solution}
\end{figure}

{Similarly to \cref{sec:accuracy}, we compare the efficiency of the three schemes on this perturbed steady solution. The reference solution is computed using the \FWBschemeTwo scheme with $2^{13}$ cells ($2^{11}$ for the \water and \methane EOSs, since they are quite costly to compute). The results are displayed in \cref{fig:efficiency_steady_solution}. This time, the \FWBscheme scheme clearly outperforms the \HLLscheme one, even on highly nonlinear and tabulated EOSs; the well-balanced correction is responsible for this gain in performance. The \FWBschemeTwo scheme remains more efficient than the \FWBscheme one, but less clearly so than in the previous case.}

\subsubsection{Sinusoidal perturbation of a moving steady solution}
\label{sec:perturbation_sin}

Another test case of interest consists in a moving steady solution perturbed by
a wave created by a time-dependent boundary condition and
propagating into the computational domain.
Similar set-ups can be found in \cite{FucMcmMisRisWaa2010},
motivated by the study of wave propagation in stellar atmospheres.
Here, we consider the same moving equilibrium as above, described by $(q_0, s_0, H_0)$ in \cref{tab:MovingEQ_Triplet}.
The perturbation is applied onto the momentum as a right (left for \PenRob) boundary condition.
It is given for $t>0$ by
\begin{equation}
    q(x_0,t) = q_0 \bigl(1 + \nu \sin(\kappa \pi t) \bigr).
\end{equation}
Therein, $\nu$ denotes the amplitude of the perturbation, $\kappa$ its frequency, and $x_0 \in \{0, 1\}$ corresponds to the excited boundary condition.
Moreover, $q_0 > 0$ is the steady momentum reported in \cref{tab:MovingEQ_Triplet} for each EOS respectively,
and so the flow travels from left to right.
Based on the same motivation as in the test cases above, these parameters are chosen differently for each EOS and summarized in \cref{tab:wavy_BC_parameters}.
The parameter $x_0$ is taken equal to $1$ to perturb the right boundary,
except for the \PenRob EOS where we take $x_0 = 0$.
Indeed, for the \PenRob EOS, all waves travel towards the right.
Therefore, perturbing the left boundary for the \PenRob EOS allows us
to observe the waves propagating through the computational domain.
For the other five EOSs, waves travel in both directions,
and thus we perturb the right boundary
to observe the waves going upwind with respect to the flow direction.
The simulation is stopped as the perturbation reaches the opposite boundary of the computational domain.
Due to the different initial conditions and configurations of the EOSs, the final time $t_\text{f}$ differs and {is} reported in \cref{tab:MovingEQ_final_times}.

\begin{table}[!bt]
    \centering
    \renewcommand{\arraystretch}{1.25}
    \begin{tabular}{rcccc}
        \toprule
        EOS      &
                 & $\nu$
                 & $\kappa \left[\unit{\hertz}\right]$
                 & $\Lambda$
        \\
        \cmidrule(lr){1-5}
        \ideal   &
                 & $10^{-8}$
                 & \num{8}
                 & \num{1.5}
        \\
        \vdW     &
                 & $10^{-8}$
                 & \num{16}
                 & \num{2.5}
        \\
        \RedKwo  &
                 & $10^{-8}$
                 & \num{16}
                 & \num{1.5}
        \\
        \PenRob  &
                 & $10^{-3}$
                 & \num{32}
                 & \num{1.5}
        \\
        \water   &
                 & $10^{-5}$
                 & \num{8e3}
                 & \num{1.1}
        \\
        \methane &
                 & \num{5e-6}
                 & \num{8e3}
                 & \num{1.1}
        \\
        \bottomrule
    \end{tabular}
    \caption{Steady solution perturbed at the boundary: required values of the parameters. The values of the perturbation constants are reported, alongside the final time and value of the wave speed factor $\Lambda$.}
    \label{tab:wavy_BC_parameters}
\end{table}

\begin{figure}[htb!]
    \centering
    \includegraphics[scale=1.25]{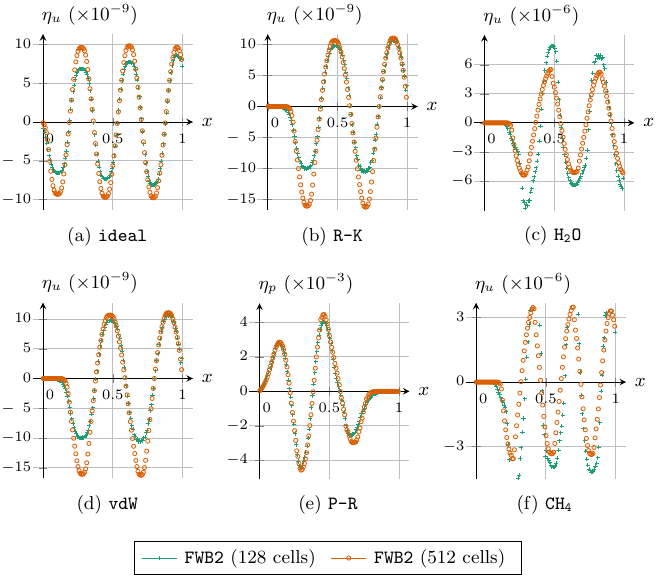}
    \caption{Boundary perturbation of equilibrium state, described in \cref{sec:perturbation_sin}: Relative velocity perturbation $\reldiff_u = (u_\text{eq}(x) - u)/u_\text{eq}(x)$ or relative pressure perturbation $\reldiff_p = (p_\text{eq}(x) - p)/p_\text{eq}(x)$, with respect to the equilibrium velocity $u_\text{eq}(x)$ or pressure $p_\text{eq}(x)$ respectively. {The results are obtained by the \FWBschemeTwo scheme with either \num{128} or \num{512} cells, using the six EOSs under consideration.}}
    \label{fig:moving_perturbed_bc}
\end{figure}

The relative difference between numerical results and the moving equilibrium
are given in \cref{fig:moving_perturbed_bc}{, for the \FWBschemeTwo scheme on two grids, made of \num{128} and \num{512} cells respectively.}
Note that we have modified the factor $\Lambda$ in the wave speed estimate, to introduce some additional numerical viscosity to stabilize the numerical scheme, see \cref{tab:wavy_BC_parameters}.
{These changes are marginal and are necessary since the
wave speed estimates \eqref{eq:def_lambda} with $\Lambda = 1$
may violate the stability conditions,
as per~\cite{TorMueSiv2020}}.

We observe that the perturbations remain well-resolved for all considered EOS despite being of orders around $10^{-9}$ in some cases.
Note that, for the \PenRob EOS, the pressure perturbation is displayed instead of the velocity perturbation.
This is due to the steady state being supercritical, and thus the perturbations in $u$ quickly leave the domain, while the acoustic waves remain visible as pressure perturbations.

\subsection{Riemann Problems}
\label{sec:RPs}

As a final series of test cases,
we consider three Riemann Problems (RPs),
two classical ones for the Euler equations without a source term,
and one in presence of a gravitational field, far from an equilibrium.

In all these cases, the space domain is $\Omega = [0,1]$,
and we prescribe homogeneous Neumann boundary conditions.
Moreover, the initial condition takes the form of a Riemann problem, i.e.,
\begin{equation*}
    W(x,0) = \begin{cases}
        W_L & \text{if } x < x_0,    \\
        W_R & \text{if } x \geq x_0,
    \end{cases}
\end{equation*}
with the jump position $x_0 = 0.5$.
The left and right states depend on the problem under consideration.
In each case, we take {a coarse grid made of $75$ cells for all schemes.
        To underline the improvements with the second-order approach, we add the solution computed with the first-order \FWBscheme on a fine grid made of $750$ cells.}

    {Furthermore, we verify the entropy stability.
        Namely, since we are in a closed domain,
        and the waves do not reach the boundaries,
        summing the discrete entropy inequality
        \eqref{discreteentropineq} over the space domain yields
        \begin{equation*}
            \sum_{i=1}^N \rho_i^{n+1}\eta(s_i^{n+1})
            \leqslant
            \sum_{i=1}^N \rho_i^n \eta(s_i^n)
            \text{,\qquad i.e., \qquad}
            \sum_{i=1}^N U_i^{n+1}
            \leqslant
            \sum_{i=1}^N U_i^n,
        \end{equation*}
        with $U = \rho \eta(s(W))$ being the mathematical entropy.
        This allows us to chart the entropy decay
        by displaying the sequence
        \begin{equation*}
            (U^n)_{n \geq 0}
            =
            \left(
            \sum_{i=1}^N U_i^{n}
            -
            \sum_{i=1}^N U_i^0
            \right)_{n \geq 0}
        \end{equation*}
        in time, which should be decreasing according to the (discrete) entropy inequality.
        We compute the mathematical entropy exemplarily by using three increasing and convex functions $\eta(s)$,
        namely,
        \begin{equation*}
            \eta_1(s) = e^s, \qquad
            \eta_2(s) = \ln(1 + e^s), \qquad
            \eta_3(s) = s + \sqrt{1 + s^2}.
        \end{equation*}
        In practice,
        we compute $\eta(s/s_\text{max})$,
        where $s_\text{max}$ is the maximum value of
        $s$ in the simulation.
        This helps to avoid scaling issues,
        especially for the tabulated EOSs,
        where the entropy is large enough for the computation of
        $e^s$ to result in an overflow.}

\subsubsection{Homogeneous case}
\label{sec:homogeneous_RP}

We first consider RPs without the influence of the gravitational field, i.e. we set $\varphi = 0$.
Note that, for the \vdW EOS,
we set $a_0 = \qty{2}{\m^6\per\kg\squared}$ and $b = \qty{0.5}{\m\cubed\per\kg}$,
to have solutions closer to the \ideal EOS,
and to check the impact of this more complex EOS on classical RPs.

\begin{table}[!bt]
    \centering
    \renewcommand{\arraystretch}{1.25}
    \begin{tabular}{cccccc}
        \toprule
        EOS                                                    &       &
        $\rho \left[\frac{\unit{\kg}}{\unit{\m\cubed}}\right]$ &
        $u \left[\frac{\unit{\m}}{\unit{\second}}\right]$      &
        $p \left[\unit{\Pa}\right]$                            &
        $t_\text{f} \, [\unit{\s}]$                                                   \\
        \cmidrule(lr){1-6}
        \multirow{2}{*}{\ideal}                                &
        L                                                      & 1     & 0 & 1      &
        \multirow{2}{*}{0.1644}                                                       \\
                                                               &
        R                                                      & 0.125 & 0 & 0.1    & \\
        \cmidrule(lr){1-6}
        \multirow{2}{*}{\vdW}                                  &
        L                                                      & 1     & 0 & 1      &
        \multirow{2}{*}{0.1644}                                                       \\
                                                               &
        R                                                      & 0.125 & 0 & 0.1    & \\
        \cmidrule(lr){1-6}
        \multirow{2}{*}{\RedKwo}                               &
        L                                                      & 1.5   & 0 & 1.25   &
        \multirow{2}{*}{0.25}                                                         \\
                                                               &
        R                                                      & 0.5   & 0 & 0.75   & \\
        \cmidrule(lr){1-6}
        \multirow{2}{*}{\PenRob}                               &
        L                                                      & 1.5   & 0 & 1.5    &
        \multirow{2}{*}{0.25}                                                         \\
                                                               &
        R                                                      & 0.35  & 0 & 0.15   & \\
        \cmidrule(lr){1-6}
        \multirow{2}{*}{\water}                                &
        L                                                      & 996   & 0 & $10^6$ &
        \multirow{2}{*}{\num{2e-4}}                                                   \\
                                                               &
        R                                                      & 998   & 0 & $10^5$ & \\
        \cmidrule(lr){1-6}
        \multirow{2}{*}{\methane}                              &
        L                                                      & 422   & 0 & $10^7$ &
        \multirow{2}{*}{\num{2e-4}}                                                   \\
                                                               &
        R                                                      & 424   & 0 & $10^5$ & \\
        \bottomrule
    \end{tabular}
    \caption{Initial conditions for the Sod-like Riemann problems. Initial states left (L) and right (R) are reported for density, velocity and pressure as well as the final time $t_\text{f}$.}
    \label{tab:RPinit}
\end{table}

\begin{figure}[!tb]
    \centering
    \centering
    \includegraphics[scale=1.25]{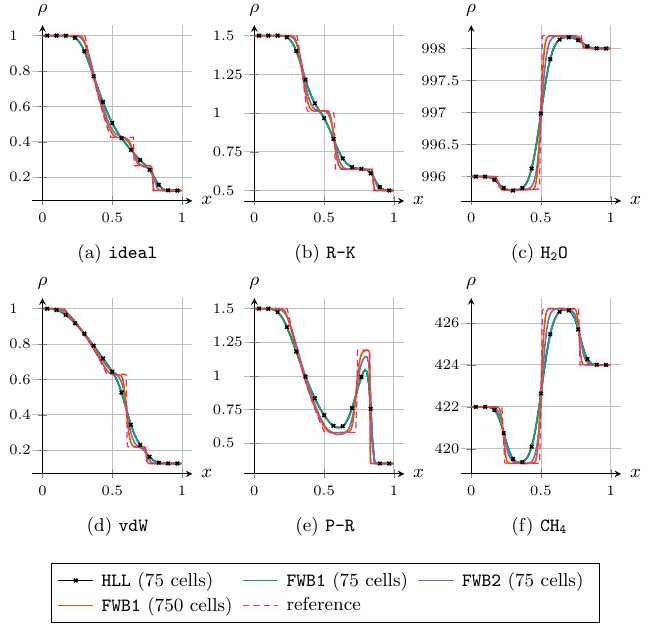}
    \caption{Sod-like Riemann Problems: Density $\rho$ at final time for the six EOSs, obtained with the \HLLscheme, \FWBscheme and \FWBschemeTwo schemes, against a reference solution. For each EOS, the final time, as well as the left and right states, are reported in \cref{tab:RPinit}.}
    \label{fig:RP_dens_no_gravity}
\end{figure}

\begin{figure}[!tb]
    \centering
    \centering
    \includegraphics[scale=1.25]{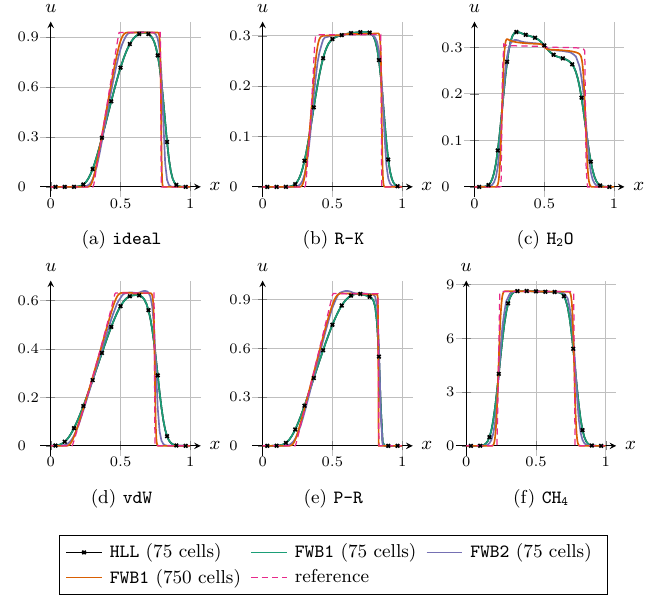}
    \caption{Sod-like Riemann Problems: Velocity $u$ at final time for the six EOSs, obtained with the \HLLscheme, \FWBscheme and \FWBschemeTwo schemes, against a reference solution. For each EOS, the final time, as well as the left and right states, are reported in \cref{tab:RPinit}.}
    \label{fig:RP_u_no_gravity}
\end{figure}

\begin{figure}[!tb]
    \centering
    \centering
    \includegraphics[scale=1.25]{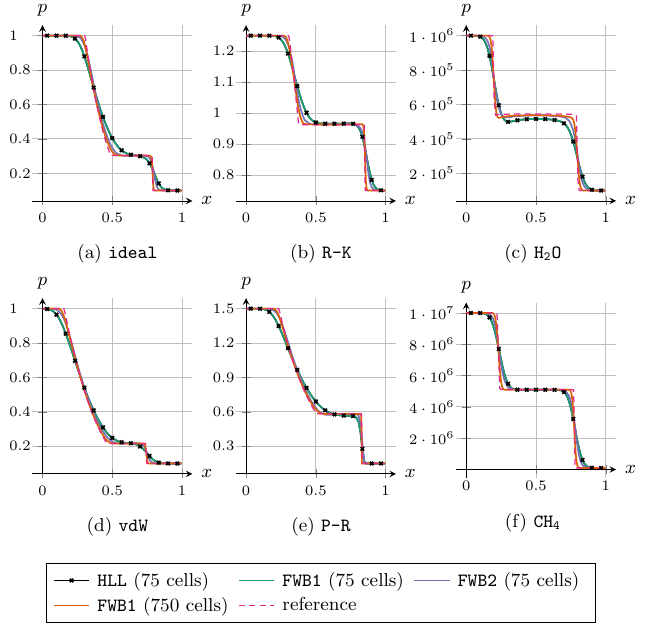}
    \caption{Sod-like Riemann Problems: Pressure $p$ at final time for the six EOSs, obtained with the \HLLscheme, \FWBscheme and \FWBschemeTwo schemes, against a reference solution. For each EOS, the final time, as well as the left and right states, are reported in \cref{tab:RPinit}.}
    \label{fig:RP_p_no_gravity}
\end{figure}

\begin{figure}[!tb]
    \centering
    \centering
    \includegraphics[scale=1.25]{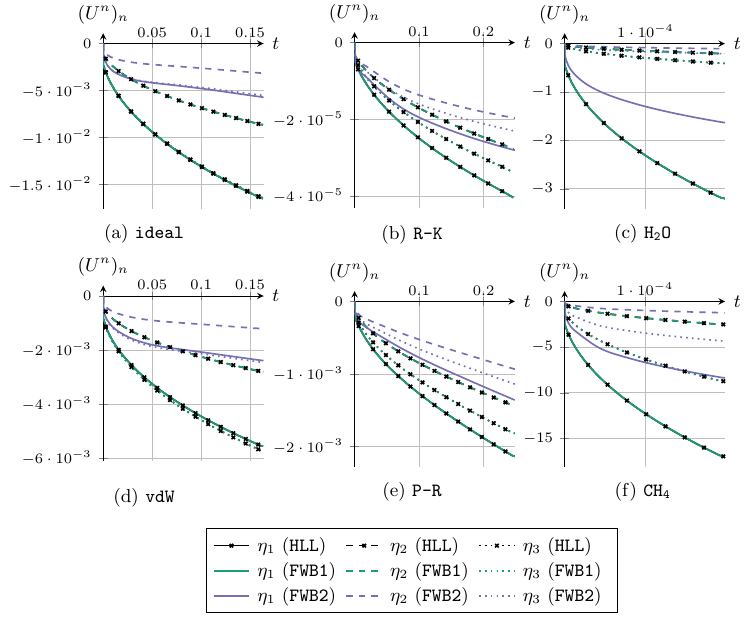}
    \caption{Sod-like Riemann Problems: Entropy decay over time for the \HLLscheme, \FWBscheme and \FWBschemeTwo schemes with $75$ cells, for the six EOSs.}
    \label{fig:RP_entropy_decay_no_gravity}
\end{figure}

We first consider a set of Sod-like RPs, whose
left and right states are given for each EOS in \cref{tab:RPinit}.
The numerical results for $\rho$, $u$ and $p$ are given in
\cref{fig:RP_dens_no_gravity,fig:RP_u_no_gravity,fig:RP_p_no_gravity}, respectively.
    {We remark that the \FWBscheme {and \FWBschemeTwo} schemes are} able to correctly determine the shock position and amplitude independently of the applied EOS.
The numerical results are compared against a reference solution obtained with the classical HLL scheme \cite{HarLaxLee1983} using \num{7500} cells.
We observe a good agreement with the reference solution.
Note, however, that the contact discontinuity has been quite
diffused in the \water case (sub-figures (c)).
This leads to a perturbation of the plateau
in the velocity profile.
In all cases, we observe that the results obtained with \HLLscheme and \FWBscheme are superimposed.
This is due to $\varphi = 0$ where both schemes are essentially the same.
    {Moreover, as expected, the \FWBschemeTwo scheme is able to capture all features of the solution more accurately than the \FWBscheme one.
        Even more, the results are comparable with the \FWBscheme scheme on the much finer grid made of $750$ cells. In \cref{fig:RP_entropy_decay_no_gravity} the entropy decay
        for the \FWBscheme, \FWBschemeTwo and \HLLscheme schemes with $75$ cells is reported.
        The total mathematical entropy $U = \rho \eta(s(W))$
        is indeed decreasing over time,
        for each of the three choices of $\eta$
        and the six EOSs.}


\begin{figure}[!hbt]
    \centering
    \centering
    \includegraphics[scale=1.25]{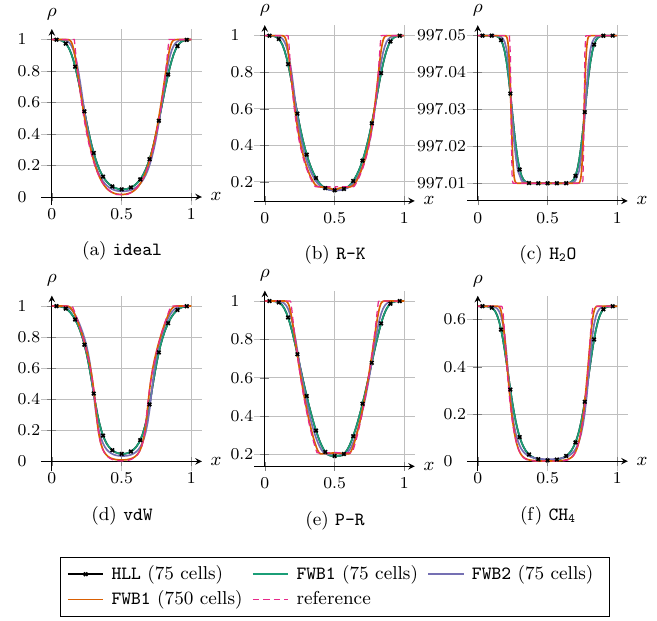}
    \caption{Double rarefaction Riemann Problems: Density $\rho$ at final time for the six EOSs, obtained with the \HLLscheme and \FWBscheme schemes, against a reference solution. For each EOS, the final time, as well as the initial condition, are reported in \cref{tab:RP_double_rarefaction_init}.}
    \label{fig:RP_double}
\end{figure}

\begin{figure}[!tb]
    \centering
    \centering
    \includegraphics[scale=1.25]{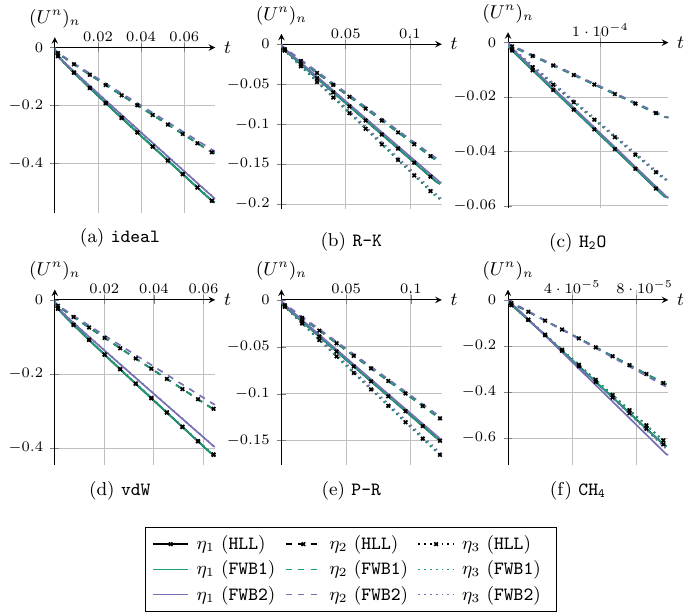}
    \caption{Double rarefaction Riemann Problems: Entropy decay over time for the \HLL, \FWBscheme and \FWBschemeTwo schemes with $75$ cells, for all six EOSs.}
    \label{fig:RP_double_entropy_decay}
\end{figure}

We then turn to a double rarefaction,
whose initial data for all EOS in terms of a Riemann problem {are} given by
\begin{equation}
    \rho_L = \rho_R = \rho_0, \quad q_L = -q_0, \quad q_R = q_0, \quad p_L = p_R = p_0.
\end{equation}
The initial data {are} detailed in \cref{tab:RP_double_rarefaction_init}.
This is a challenging problem,
as it can lead to a near-vacuum state in the center of the domain.
The results are displayed for the density in \cref{fig:RP_double}
at the final times reported in \cref{tab:RP_double_rarefaction_init}.
A similar behavior is obtained for the pressure.
However, due to the similarity in the results, we have omitted the pressure plots.
Despite the fact that $\rho$ and~$p$ are very close to zero in the center of the domain,
see e.g. the \ideal, \vdW or \methane EOS,
no negative values of $\rho$ and $p$ are observed nor the simulation had to be stopped.
The numerical results are compared against a reference solution obtained with the classical \HLLscheme solver \cite{HarLaxLee1983} using \num{7500} cells.
The entropy decay
for the \FWBscheme, {\FWBschemeTwo and \HLL} schemes with $75$ cells
are displayed in \cref{fig:RP_double_entropy_decay}.
We indeed observe that the discrete entropy inequality
is satisfied for all choices of $\eta$ and all EOSs.

\begin{table}[!hbt]
    \centering
    \renewcommand{\arraystretch}{1.25}
    \begin{tabular}{rccccc}
        \toprule
        EOS      &
                 & $\rho_0 \left[\frac{\unit{\kg}}{\unit{\m\cubed}}\right]$
                 & $q_0 \left[\frac{\unit{\kg}}{\unit{\square\m\second}}\right]$
                 & $p_0 \left[\unit{\Pa}\right]$
                 & $t_\text{f} \, [\unit{\s}]$                                   \\
        \cmidrule(lr){1-6}
        \ideal   &
                 & 1
                 & 10/3
                 & 1
                 & 0.075                                                         \\
        \vdW     &
                 & 1
                 & 10/3
                 & 1
                 & 0.065                                                         \\
        \RedKwo  &
                 & 1
                 & 1.8
                 & 1
                 & 0.125                                                         \\
        \PenRob  &
                 & 1
                 & 1.5
                 & 1
                 & 0.125                                                         \\
        \water   &
                 & 997.05
                 & 59.9
                 & \num{101325}
                 & \num{1.75e-4}                                                 \\
        \methane &
                 & 0.657
                 & \num{1789}
                 & \num{101325}
                 & $10^{-4}$                                                     \\
        \bottomrule
    \end{tabular}
    \caption{Initial conditions for the double rarefaction Riemann problems. The initial density, momentum and pressure is reported, as well as the final time $t_\text{f}$, for each EOS.}
    \label{tab:RP_double_rarefaction_init}
\end{table}


\subsubsection{Riemann problem in a gravitational field}

\begin{figure}[!tb]
    \centering
    \centering
    \includegraphics[scale=1.25]{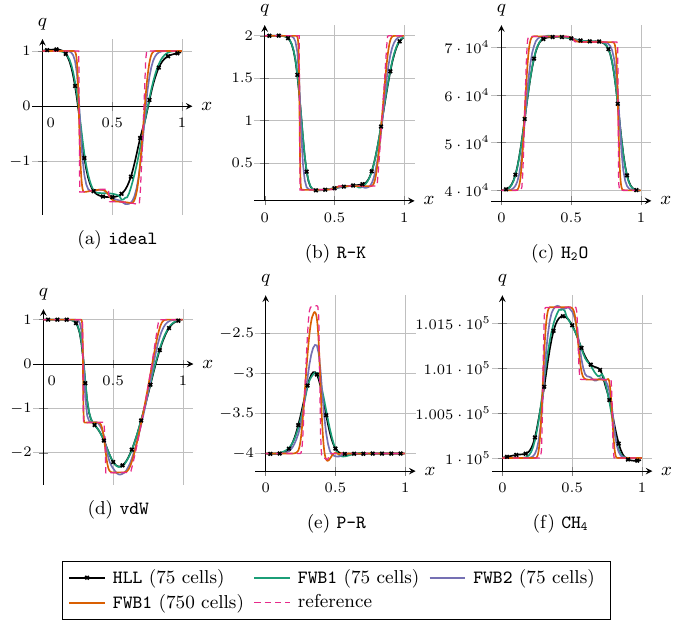}
    \caption{Sod-like Riemann Problems in a gravitational field: Momentum $q$ at final time for the six EOSs, obtained with the \HLLscheme and \FWBscheme schemes, against a reference solution. For each EOS, the final time, as well as the left and right states, are reported in \cref{tab:RP_source_init}.}
    \label{fig:RP_dens}
\end{figure}

\begin{figure}[!tb]
    \centering
    \centering
    \includegraphics[scale=1.25]{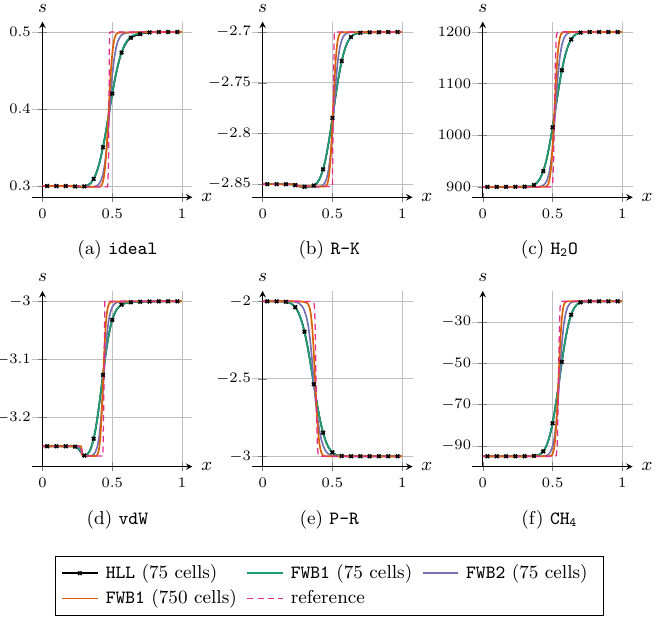}
    \caption{Sod-like Riemann Problems in a gravitational field: Specific entropy $s$ at final time for the six EOSs, obtained with the \HLLscheme and \FWBscheme schemes, against a reference solution. For each EOS, the final time, as well as the left and right states, are reported in \cref{tab:RP_source_init}.}
    \label{fig:RP_u}
\end{figure}

\begin{figure}[!tb]
    \centering
    \centering
    \includegraphics[scale=1.25]{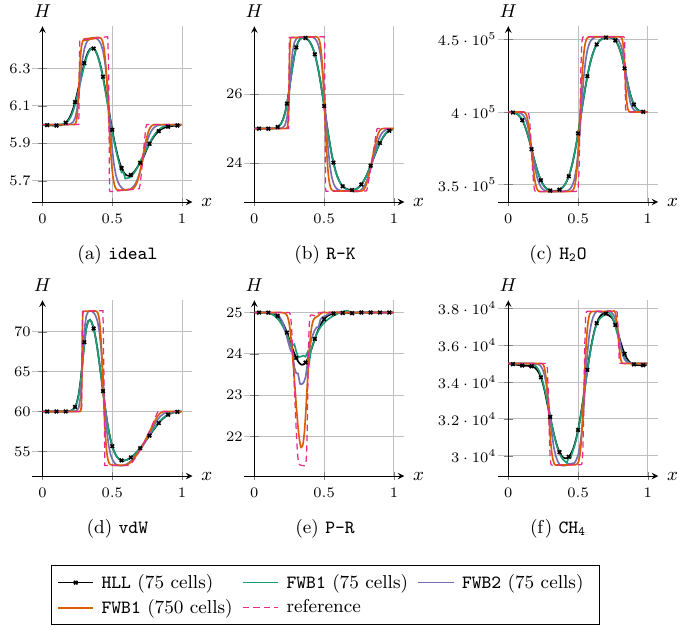}
    \caption{Sod-like Riemann Problems in a gravitational field: Specific total enthalpy $H$ at final time for the six EOSs, obtained with the \HLLscheme and \FWBscheme schemes, against a reference solution. For each EOS, the final time, as well as the left and right states, are reported in \cref{tab:RP_source_init}.}
    \label{fig:RP_p}
\end{figure}

\begin{figure}[!tb]
    \centering
    \centering
    \includegraphics[scale=1.25]{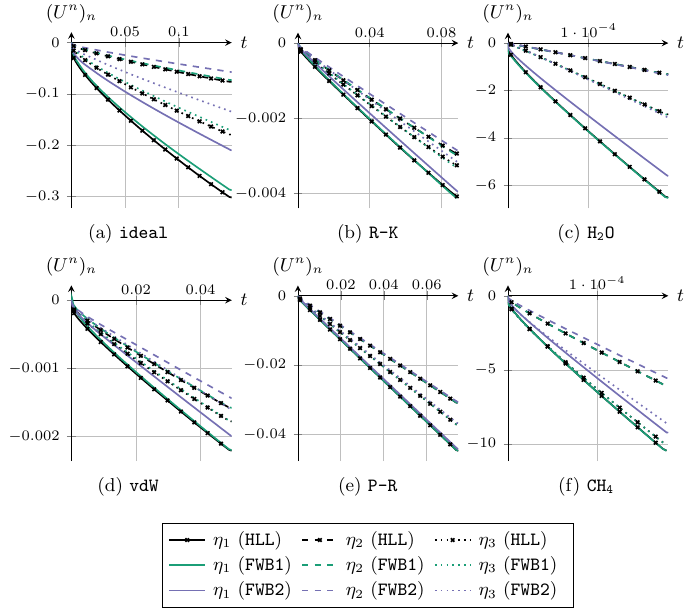}
    \caption{Sod-like Riemann Problems in a gravitational field: Entropy decay over time for the \HLLscheme, \FWBscheme and \FWBschemeTwo schemes with $75$ cells, for the six EOSs.}
    \label{fig:RP_entropy_decay}
\end{figure}

{The last Riemann problem is a Sod-like problem in the presence of the gravitational potential given by
    \begin{equation*}
        \varphi(x) =
        \begin{dcases}
            \frac{\varphi_0}{4}
            \exp \left( 1 - \frac 1 {1 - (4x - 2)^2} \right)
             & \text{if } \left|x - \frac 1 2\right| < \frac 1 4, \\
            0
             & \text{otherwise.}
        \end{dcases}
    \end{equation*}
    This compactly supported potential is used to
    make sure that the domain is closed,
    to observe the entropy decay.}
The initial data {are} given in steady variables $q$, $s$ and $H$,
and {are} reported in \cref{tab:RP_source_init},
alongside the final time $t_\text{f}$, for each EOS.
Note that, for the \PenRob EOS, we set the initial jump position
to $x_0 = 0.8$, since all waves travel
towards the left of the domain.
The equilibrium variables are displayed in \cref{fig:RP_dens,fig:RP_u,fig:RP_p}.
We observe that the solution generated by the new well-balanced \FWBscheme {and \FWBschemeTwo schemes} accurately captures the arising shock and rarefaction waves.
    {The entropy decay, displayed in \cref{fig:RP_entropy_decay}, shows that the schemes also respect the discrete entropy inequality.}
The numerical results are compared against a reference solution
obtained with the \HLLscheme scheme on \num{7500} cells.
    {We also compare our results to the
        \HLLscheme scheme on a coarse mesh.
        This time, the curves are not superimposed,
        due to the difference in source term treatment.
        Most notably, we observe that the \FWBscheme and \FWBschemeTwo
        are able to capture plateaux and waves
        in a more accurate way than the \HLLscheme scheme.
        This is confirmed by comparing the \FWBscheme scheme on a coarse mesh (green curves) to its counterpart on a fine mesh (orange curves).
        In the latter case, the plateaux are sharper and the waves are better resolved.}

\begin{table}[bt!]
    \centering
    \newcommand{\unitQ}{\left[\frac{\unit{\kg}}{\unit{\square\m\second}}\right]}
    \newcommand{\unitS}{\left[\frac{\unit{\joule}}{\unit{\kelvin}}\right]}
    \newcommand{\unitH}{\left[\unit{\joule}\right]}
    \renewcommand{\arraystretch}{1.25}
    \begin{tabular}{cccccc}
        \toprule
        EOS                       &                            &
        $q \unitQ$                &
        $s \unitS$                &
        $H \unitH$                &
        $t_\text{f} \, [\unit{\s}]$                                                                     \\
        \cmidrule(lr){1-6}
        \multirow{2}{*}{\ideal}   &
        L                         & \multirow{2}{*}{1}         & 0.3   & \multirow{2}{*}{6}           &
        \multirow{2}{*}{0.15}                                                                           \\
                                  &
        R                         &                            & 0.5   &                              & \\
        \cmidrule(lr){1-6}
        \multirow{2}{*}{\vdW}     &
        L                         & \multirow{2}{*}{1}         & -3.25 & \multirow{2}{*}{60}          &
        \multirow{2}{*}{0.05}                                                                           \\
                                  &
        R                         &                            & -3    &                              & \\
        \cmidrule(lr){1-6}
        \multirow{2}{*}{\RedKwo}  &
        L                         & \multirow{2}{*}{2}         & -2.85 & \multirow{2}{*}{25}          &
        \multirow{2}{*}{0.09}                                                                           \\
                                  &
        R                         &                            & -2.7  &                              & \\
        \cmidrule(lr){1-6}
        \multirow{2}{*}{\PenRob}  &
        L                         & \multirow{2}{*}{-4}        & -2    & \multirow{2}{*}{25}          &
        \multirow{2}{*}{0.075}                                                                          \\
                                  &
        R                         &                            & -3    &                              & \\
        \cmidrule(lr){1-6}
        \multirow{2}{*}{\water}   &
        L                         & \multirow{2}{*}{\num{4e4}} & 900   & \multirow{2}{*}{\num{4e5}}   &
        \multirow{2}{*}{\num{2e-4}}                                                                     \\
                                  &
        R                         &                            & 1200  &                              & \\
        \cmidrule(lr){1-6}
        \multirow{2}{*}{\methane} &
        L                         & \multirow{2}{*}{\num{1e5}} & -95   & \multirow{2}{*}{\num{3.5e4}} &
        \multirow{2}{*}{\num{1.8e-4}}                                                                   \\
                                  &
        R                         &                            & -20   &                              & \\
        \bottomrule
    \end{tabular}
    \caption{Initial conditions for the Sod-like Riemann problems in a gravitational field. Initial states left (L) and right (R) are reported for momentum $q$, specific entropy $s$ and specific total enthalpy $H$; the final time $t_\text{f}$ is also reported.}
    \label{tab:RP_source_init}
\end{table}


\section{Conclusions}
\label{sec:Conclusion}
In this paper, we derived a new numerical scheme to
approximate weak solutions of the Euler equations with a gravitational source term with general equations of state which can be given analytically or tabulated.
The presented Godunov-type finite volume scheme is based on an approximate Riemann solver composed of two intermediate states.
It is constructed such that the scheme is consistent, fully well-balanced and  preserves the positivity of the density.
Moreover, if the system is equipped with a convex entropy and an associated entropy inequality holds, the scheme is entropy-stable and the positivity of thermodynamic variables is ensured.
The proof of these properties are summarized in \cref{theo:summary}.
The results presented in this work are an extension of the fully well-balanced scheme for the ideal gas law detailed in Berthon et al. \cite{Berthon2025}.
The theoretical findings are validated numerically by computing moving equilibrium solutions as well as the propagation of perturbations and Riemann Problems using three different analytically given equations of states as well as two tabulated equations of state based on the \texttt{CoolProp} library \cite{CoolProp}.

\section*{Acknowledgments}

V. Michel-Dansac extends his thanks to ANR-22-CE25-0017 OptiTrust.
The SHARK-FV conference has greatly contributed to this work.
    {We would like to thank the anonymous referees whose remarks and comments have greatly contributed to the quality of the manuscript.}

\end{document}